\documentclass[11pt]{article}

\usepackage{epsfig}
\usepackage{amsmath}
\usepackage{soul,color}
\usepackage{bm}
\usepackage{epsfig}
\usepackage{natbib}
\newcommand{\tr}{^{\prime}}
\def\b#1{\mbox{\boldmath $#1$}}    

\renewcommand{\th}{\theta}

\newcommand{\al}{\alpha}
\newcommand{\be}{\beta}

\newcommand{\si}{\sigma}

\newcommand{\ep}{\varepsilon}
\def\bbeta{\mbox{\boldmath$\beta$}}

\def\bSigma{\mbox{\boldmath$\Sigma$}}

\usepackage{bm}
\usepackage{epsfig}

\begin{document}
\title{\vspace*{-1cm} Mixture latent autoregressive models for longitudinal data}
\author{Francesco Bartolucci\footnote{Department of Economics, Finance and
Statistics, University of Perugia, Via A. Pascoli, 20, 06123
Perugia, Italy.} \footnote{{\em e-mail}: bart@stat.unipg.it} ,\hspace{1mm}
Silvia Bacci$^*$\footnote{{\em e-mail}: sbacci@stat.unipg.it} ,\hspace{1mm}
Fulvia Pennoni\footnote{Department of Statistics,
University of Milano-Bicocca, Via Bicocca degli Arcimboldi 8, 20126 Milano, Italy.}
\footnote{{\em e-mail}: fulvia.pennoni@unimib.it}}
\date{}\maketitle
%
\begin{abstract}
\noindent Many relevant statistical and econometric models for the
analysis of longitudinal data include a latent process to account
for the unobserved heterogeneity between subjects in a dynamic
fashion. Such a process may be continuous (typically an AR(1)) or
discrete (typically a Markov chain). In this paper, we propose a
model for longitudinal data which is based on a mixture of AR(1)
processes with different means and correlation coefficients, but
with equal variances. This model belongs to the class of models
based on a continuous latent process, and then it has a natural
interpretation in many contexts of application, but it is more
flexible than other models in this class, reaching a goodness-of-fit
similar to that of a discrete latent process model, with a reduced
number of parameters. We show how to perform maximum likelihood
estimation of the proposed model by the joint use of an
Expectation-Maximisation algorithm and a Newton-Raphson algorithm,
implemented by means of recursions developed in the hidden Markov
literature. We also introduce a simple method to obtain standard
errors for the parameter estimates and a criterion to choose the
number of mixture components. The proposed approach is illustrated
by an application to a longitudinal dataset, coming from the Health
and Retirement Study, about self-evaluation of the health status by
a sample of subjects. In this application, the response variable is
ordinal and time-constant and time-varying individual covariates are
available.

\end{abstract}

\newpage

\section{Introduction}\label{sec:Introduction}
In the analysis of longitudinal data, an important aspect that can
be accounted for is the unobservable heterogeneity between subjects.
This form of heterogeneity corresponds to the effect that
unobservable factors have on the occasion-specific response
variables in addition to the effect of observable covariates. The
simplest approach to account for the unobserved heterogeneity is
based on the inclusion, in the model of interest, of
individual-specific random intercepts, that can have 
either a continuous or a discrete distribution. Models based on individual
parameters having a continuous distribution may be casted in the
class of Generalised 
Linear Mixed models and that of 
Random Effects models \citep[see][]{sni:bos:99,McCu:Sear:01,gol:03,
skr:rabe:04,hanc:samu:08}. Models based on discrete random effects
may be seen as forms of Latent Class (LC) models \citep{laza:50,
laza:henr:68, good:74,band:migl:zege:rath:97,huan:band:04}.
See also \cite{hage:mccu:02} for an
exhaustive review about the LC model.

The approaches mentioned above assume that the effect of
unobservable factors on the response variables is time constant. A
more general assumption consists of introducing, for each subject,
time-varying individual random effects which give rise to a latent
process for the unobserved heterogeneity. Even in this case we can
disentangle the continuous case from the discrete case. The most
common formulation based on a continuous-valued latent process
assumes that the individual effects follow an Autoregressive model
of order 1 (AR(1)); see \cite{chi:89} and \cite{hei:08}. Hereafter,
this model is referred to as Latent Autoregressive (LAR) model. On
the other hand, models based on a discrete latent process typically
assume that the individual effects follow a first-order Markov
chain. A Latent Markov (LM) model \citep{wigg:73} with covariates
results; see \cite{bart:farc:penn:10} for a review.

The debate on which is more appropriate, between the continuous and
the discrete latent process formulation, is open. The first formulation
is usually more easy to justify from a theoretical point of view; in
principle, there is no reason to consider the effect of unobserved
factors as discrete. Moreover, a LAR model has a parsimony close to
that of the corresponding continuous random effect model (with
time-constant individual effects), since it represents the latent
structure by only two parameters: correlation index and variance of
the individual effects. However, the model estimation may be
computationally problematic \citep{hei:08}. On the other hand,
discrete latent variable models may reach a better fit to the
analysed data. In particular, the LM model may be seen as a
semi-parametric model because, with the suitable number of states,
the underlying Markov chain may approximate any (even continuous)
process 
with a first order dependence structure. 
This advantage is at the cost of a reduced parsimony, since
the number of parameters increases with the square of the number of
states. Moreover, the interpretation may be more difficult for the
same reason mentioned above: it is more natural to consider the
effect of unobservable factors or covariates as continuous than
discrete.

A debate similar to that described above, between a continuous and a
discrete formulation for the latent process, is also present in the
literature on models for item responses 
\citep{ham:swam:85}; see, for among others, \cite{lin:91}. A similar
debate is also present in the literature about the analysis of
certain types of time-series data. In particular, for the analysis
of financial data, the Stochastic Volatility (SV) model
\citep{taylor:1982,shep:96} may be used as an alternative to the
hidden Markov \citep[HM,][]{mac:zuc:97} model. For an interesting
comparison between the two approaches see \cite{tayl:99} and for a
comprehensive review see \cite{tayl:05}. The SV model represents the
volatility by an AR(1) process, and then has a structure that
recalls that of the LAR model for longitudinal data, whereas the HM
model relies on a Markov chain, and then it 
is very similar to the LM
model. In the field of time-series data, which is strongly related
to that of longitudinal data, we also have to mention
Markov-switching models \citep[see][Ch.~9-10]{Fruh:06}.
In this field, the interest is also on the the possibility to
combine a continuous with a discrete approach; see, among others,
\cite{kita:87}, \cite{cai:94}, \cite{ham:sus:94},
\cite{so:lam:li:98}, and \cite{ros:gal:06}. Although the
similarities between the typical formulations of models for time
series and longitudinal data, to our knowledge no attempts to
combine continuous and discrete process formulations have been made
in the context of longitudinal data.\newpage

In this paper, we propose a model for longitudinal data which is based
on a mixture of latent 
AR(1) processes to account for the unobserved heterogeneity
between subjects in a dynamic fashion. Each component
of the mixture has its
own mean and correlation coefficient, but these components have a
common variance. The proposed model, which can be used with response
variables of a different nature (binary, ordinal, or continuous),
belongs to the class of continuous latent process models for
longitudinal data. As such, it retains the natural interpretation
that characterises the LAR model, but it reaches a better fit to the
data, since it generalises this model. In particular, the goodness
of fit to the data may reach levels close that those of the LM
model, but with a reduced number of parameters.

In order to make inference on the proposed model, we show how to
compute its likelihood function by a recursion taken for the
hidden Markov literature \citep{bau:70,mac:zuc:97}; 
a procedure results which is equivalent to the
Sequential Gaussian Quadrature (SGQ) method proposed by \cite{hei:08}.
Through recursions similar to those used for HM models, we also implement an
Expectation-Maximisation (EM) algorithm and a Newton-Raphson (NR)
algorithm for the maximisation of the likelihood function and,
therefore, for the 
estimation of the model parameters.
The NR algorithm is based on the observed information matrix
which is obtained by the same numerical method proposed by \cite{bart:farc:09}.
This matrix is also exploited to obtain standard errors for the
parameter estimates.
Finally, we show how to obtain the
prediction of the individual effect for every subject in the sample
at each time occasion. Through these predictions we define a
criterion to choose the number of components to be used in 
data analysis. We recall that each component corresponds to a
separate AR(1) latent process.

The advantages of 
the proposed approach are illustrated through an
application to a longitudinal dataset concerning the self-evaluation
of health status at eight different time periods. The dataset is
derived from the the Health and Retirement Study conducted by the
University of Michigan. In this case, the response variable observed
at each occasion has five ordered categories. The proposed model is
implemented by specifying a proportional odds model for global
logits. Some observed covariates related with individual
characteristics are also 
included. The model selected for these data is
compared with the corresponding LAR and LM models.\newpage

The paper is organised as follows. In the next section we introduce
the basic notation and describe some relevant approaches for
longitudinal data. In Section 3 we outline the proposed model for
longitudinal data and, in Section 4, we describe likelihood based
inference for this model. The results of the application based on
the self evaluation of health status data are illustrated in Section
5. Final conclusions are reported in the Section 6.
\section{Preliminaries}\label{sec2}
With reference to a sample of $n$ subjects observed at $T$ time
occasions, let $y_{it}$ be the response variable for subject
$i$ at occasion $t$ and let $\b x_{it}$ be a
corresponding column 
vector of covariates, with $i=1,\ldots,n$ and $t=1,\ldots,T$.
We also denote by $\b
y_i=(y_{i1},\ldots,y_{iT})$ the vector of response variables and by
$\b X_i=\begin{pmatrix}\b x_{i1} & \cdots & \b x_{iT}\end{pmatrix}$
the matrix of all covariates for subject $i$.

At this stage, we do not restrict the response variable to have a
specific nature. Therefore, we introduce a latent continuous
variable $y_{it}^*$ underlying each $y_{it}$. In particular, we
assume that
\begin{equation}\label{eq:link}
y_{it}=G(y_{it}^*),
\end{equation}
where $G(\cdot)$ is a parametric function  which may depend on
specific parameters according to the different nature of $y_{it}$,
such as specific cutpoints in the presence of ordinal variables. 
Hereafter, we use $\b y_i^*=(y_{i1}^*,\ldots,y_{iT}^*)$ to
denote the vector of latent response variables corresponding to $\b
y_i$. More details about the possible formulations of $G(\cdot)$
will be given in Section \ref{sec:parameterizations}.

In the following, we briefly review models which allow us to
take into account the unobserved heterogeneity between subjects by introducing
time-constant and time-varying individual effects.
\subsection{Models with time-constant individual effects}
The simplest approach to take into account the unobserved
heterogeneity is based on the assumption that, for every unit 
$i$, the latent
response variables 
in $\b y_i^*$ are conditionally independent given  the
covariates $\b X_i$ and an individual-specific intercept $\al_i$
({\em local independence}).
Moreover, it is assumed that each $y_{it}^*$ only depends on $\al_i$
and $\b x_{it}$ as follows
\begin{equation}  \label{eq:randint}
y_{it}^* = \al_i+\b x_{it}\tr\bbeta +\ep_{it},\quad i=1,\ldots,n,\:
t=1,\ldots,T,
\end{equation}
where the error terms $\ep_{it}$
are assumed mutually independent and identically distributed. Note that each vector
$\b x_{it}$ may also include the lagged response.

When the individual-specific intercepts $\al_i$ are treated as
random parameters, the same distribution (usually independent of the
covariates) is assumed for all subjects, which may be continuous or
discrete. In the first case, we typically assume that $\al_i\sim
N(0,\si^2)$ for $i=1,\ldots,n$. In the second case, instead, every
$\al_i$ may assume a value among $k$ possible values or support
points $\xi_h$ having probabilities $\pi_h$, with $h=1,\ldots,k$.
The support points and the corresponding probabilities are typically
estimated on the basis of the data, which also drive the choice of
$k$. A first model of this type is known as finite mixture of
regression models, which is an extension of a mixture of normal
distributions with averages expressed as functions of the
explanatory variables \citep{Quan:72, Quan:Rams:78}.

In any case, the assumption of conditional independence of the response variables
given the individual-specific intercepts and the covariates allows us to write
\begin{equation}
p(\b y_i|\al_i,\b X_i)=\prod_t p(y_{it}|\al_i,\b x_{it}),\quad i=1,\ldots,n,
\label{eq:conditional_distibution}
\end{equation}
where  $p(y_{it}|\al_i,\b x_{it})$ denotes the probability mass or density function of
$y_{it}$, given $\al_i$ and $\b x_{it}$, which, in turn, depends on the adopted
parameterisation; 
see equations (\ref{eq:link}) and (\ref{eq:randint}).
Then, under the random effect approach, the
{\em manifest distribution} of $\b y_i$ given $\b X_i$ is obtained by marginalising
the probability or density in (\ref{eq:conditional_distibution}) 
with respect to $\al_i$. With continuous random effects, we have
\begin{equation*}
p(\b y_i|\b X_i)=\int p(\b y_i|\al_i,\b X_i)f(\al_i)d\al_i=
\int\bigg[\prod_t p(y_{it}|\al_i,\b x_{it})\bigg]f(\al_i)d\al_i,
\label{eq:marginal_distibution_continuous}
\end{equation*}
where $f(\al_i)$ is the probability density function of every 
$\al_i$, which may depend on a specific parameter vector. With discrete random
effects, instead, we have
\begin{equation*}
p(\b y_i|\b X_i)=\sum_h\bigg[\prod_t p(y_{it}|\xi_h,\b x_{it})\bigg]\pi_h.
\label{eq:marginal_distibution_discrete}
\end{equation*}
This distribution is the base for constructing a marginal likelihood to be maximised
in order to estimate the model parameters.
\subsection{Models based on time-varying individual effects}\label{sec:time-varying}
The main drawback of the individual-specific random
intercept models described above is that they assume the effect of unobservable
factors to be time constant. This assumption may be relaxed by the
inclusion of individual-time-specific effects
$\al_{it}$, $i=1,\ldots,n$, $t=1,\ldots,T$.
In an obvious way, the assumption of local independence is extended by assuming
that, for all sample units $i$, the latent response variables in $\b y_i^*$ are
conditionally independent given $\b\al_i=(\al_{i1},\ldots,\al_{iT})$ and $\b X_i$.
Moreover, assumption (\ref{eq:randint}) is naturally extended as follows
\begin{equation}
y_{it}^*=\al_{it}+\b x_{it}\tr\bbeta +\ep_{it},\quad i=1,\ldots,n,\:
t=1,\ldots,T.\label{eq:assumption1}
\end{equation}

Given the model complexity, time-varying effects may only be assumed to be random
and not fixed. Again, two alternative approaches, continuous and  discrete, are available in
the literature. The most common continuous random-effects approach assumes that every
$\b\al_i$ follows an AR(1) process, so that
\begin{eqnarray*}
\al_{i1}&=&\ep_{i1},\label{eq:ar11}\\
\al_{it}&=&\al_{i,t-1}\rho+\ep_{it}\sqrt{1-\rho^2},
\quad t=2,\ldots,T,\label{eq:ar12}
\end{eqnarray*}
where $\ep_{it}\sim N(0,\sigma^2)$, $t = 1, \ldots, T$. This is the LAR
formulation already mentioned in Section \ref{sec:Introduction},
which was studied in detail by \cite{hei:08}.

The discrete latent process formulation assumes that, for all $i$,
$\b\al_i$ follows a 
first-order homogenous Markov chain with $k$ states denoted by $\xi_1,\ldots,\xi_k$.
This chain has initial probabilities $\pi_h$ and
transition probabilities $\pi_{h_1h_2}$,
with
\begin{eqnarray}
\pi_h&=&p(\al_{i1}=\xi_h),\quad h=1,\ldots,k,\label{eq:prob_init}\\
\pi_{h_1h_2}&=&p(\al_{i,t-1}=\xi_{h_1},
\al_{it}=\xi_{h_2}),\quad h_1,h_2=1,\ldots,k,\: t=2,\ldots,T.
\label{eq:prob_trans}
\end{eqnarray}
In other words, it is assumed that every $\al_{it}$ is conditionally
independent of $\al_{i1},\ldots, \al_{i, t-2}$ given $\al_{i,t-1}$,
but apart from this assumption, the distribution of $\b\al_i$ is
unconstrained. On the other hand, this greater flexibility
corresponds to a higher number of parameters to estimate with
already $k\geq 2$. In fact, the number of parameters involved in the
Markov chain (support points and initial and transition
probabilities) is equal to 
$(k-1)+(k-1)+k(k-1)=k^2+k-2$, taking into account the
constraints $\sum_h\pi_h=1$ and $\sum_{h_2}\pi_{h_1h_2}=1$,
$h_1=1,\ldots,k$, and that to ensure identifiability one constraint has to
be put on the support points. This is a formulation of LM type, which was
exploited by \cite{bart:farc:09} to propose a flexible class of
models for multivariate categorical longitudinal data. We have to
mention that, in order to make easier the comparison between the LAR
and the LM model, we can require that the initial probabilities
$\pi_h$ in (\ref{eq:prob_init}) coincide with those of the
stationary distribution of the chain. In this case, we have a
moderate reduction of the number of parameters which becomes equal
to $k^2-1$. 

Under both the continuous and the discrete latent process formulations,
the assumption of local independence implies that
\[
p(\b y_i|\b\al_i,\b X_i)=\prod_t p(y_{it}|\al_{it},\b x_{it}).
\]
Moreover, under the first formulation, which leads to the LAR model,
the manifest distribution of $\b y_i$ given $\b X_i$ has probability mass (or density)
function
\begin{equation}
p(\b y_i|\b X_i)=\int p(\b y_i|\b\al_i,\b X_i)f(\b\al_i)d\b\al_i.
\label{eq:marginal_distibution_continuous2}
\end{equation}
This is an integral over the $T$-dimensional space of $\b\al_i$ that
may be difficult to compute in practice. At this aim,
we can use the SGQ method proposed by Heiss (2008), which is
essentially based on rewriting expression
(\ref{eq:marginal_distibution_continuous2}) as follows
\begin{eqnarray}
p(\b y_i|\b X_i)&=&\int p(y_{i1}|\al_{i1},\b x_{i1})f(\al_{i1})
\int p(y_{i2}|\al_{i2},\b x_{i2})f(\al_{i2}|\al_{i1})\cdots\nonumber\\
&&\cdots\int p(y_{iT}|\al_{iT},\b x_{iT})f(\al_{iT}|\al_{i,T-1})
d\al_{iT}\cdots d\al_{i2}d\al_{i1}\label{eq:mani_lar}
\end{eqnarray}
and then  sequentially computing the integral involving each single
random effect $\al_{it}$, where $f(\al_{i1})$ refers to the
distribution of $\al_{i1}$ and $f(\al_{it}|\al_{i,t-1})$ to the
distribution of $\al_{it}$ given $\al_{i,t-1}$, $t=2,\ldots,T$.

When a latent Markov chain is assumed, the manifest distribution of $\b y_i$ given
$\b X_i$ is defined as follows
\begin{equation}
p(\b y_i|\b X_i)=\sum_{h_1} p(y_{i1}|\xi_{h_1},\b x_{i1})\pi_{h_1}
\sum_{h_2} p(y_{i2}|\xi_{h_2},\b x_{i2})\pi_{h_1h_2}\cdots
\sum_{h_T}p(y_{iT}|\xi_{h_T},\b x_{iT})\pi_{h_{T-1}h_T}.\label{eq:manifest_LM}
\end{equation}
In order to efficiently compute this sum, we can exploit a forward
recursion \citep{bau:70,demp:lair:rubi:77} which is well known in
the hidden Markov literature \citep{mac:zuc:97}. See \cite{bart:06}
and \cite{bart:farc:penn:10} for an efficient implementation in
matrix notation.
\section{Proposed model}\label{sec3}
In this section, we describe the proposed model for longitudinal
data, which is based on a mixture of AR(1) processes to account for
the unobserved heterogeneity in a dynamic fashion. We name this
model as Mixture Latent Autoregressive model, indicated for short by
MLAR or by MLAR($k$) when we want to mean a specific number of
mixture components $k$.
\subsection{Model assumptions}\label{sec:proposed_approach}
The proposed model is based on the following assumptions for
$i=1,\ldots,n$:
\begin{itemize}
\item[A1:] the latent response variables in $\b y_i^*$, and therefore
the observed response variables in $\b y_i$, are conditionally independent given
$\b X_i$ and a latent process $\b\al_i=(\al_{i1},\ldots,\al_{iT})$;
\item[A2:] every response variable $y_{it}^*$ in $\b y_i^*$, and then every
$y_{it}$ in $\b y_i$, only depends on $\al_{it}$ and $\b x_{it}$ through
a parameterisation formulated 
on the basis of (\ref{eq:link}) and (\ref{eq:assumption1});
\item[A3:] the latent process $\b\al_i$ has distribution given by a mixture of $k$
AR(1) processes with common variance $\si^2$.
\end{itemize}

Assumptions A1 is the usual assumption of local independence already
discussed in Section \ref{sec:time-varying}; the other two
assumptions are discussed in detail below.
\subsubsection{Assumption A2}\label{sec:parameterizations}
The introduction of an underlying continuous outcome $y_{it}^*$
related to the observed response variable $y_{it}$ as specified in
(\ref{eq:link}), allows us to adapt the model 
to several situations. Indeed, depending on the assumed distribution
for the errors $\ep_{it}$ in (\ref{eq:assumption1}) 
and on the specification of $G(\cdot)$
different models result.

The simplest case is when we let $G(y^*)=y^*$, that is the identity function, and
$\ep_{it}\sim N(0,\si^2)$ for all $i$ and $t$. In this case, a model results in which
\[
y_{it}|\al_{it},\b x_{it} \sim
N(\al_{it}+\b x_{it}\tr\bbeta,\si^2),\quad i=1,\ldots,n,\: t=1,\ldots,T.
\]
This is the typical formulation adopted with continuous response variables.

When the response variables are binary, that is $y_{it}=0,1$, we
typically assume that $G(y^*)= I\{y^*>0\}$, where $I\{\cdot\}$ is an
indicator function assuming value 1
when its argument is true and value 0 otherwise. Depending on 
the distribution of the error term $\ep_{it}$ in  model (\ref{eq:assumption1}),
a logit or probit parameterisation results. More precisely,
if we assume a logistic distribution for the error terms $\ep_{it}$, then
a logit parameterisation results, under which
\begin{equation}\label{eq:logit}
\log \frac{p(y_{it}=1|\al_{it},\b x_{it})}{p(y_{it}=0|\al_{it},\b
x_{it})} = \al_{it}+\b x_{it}\tr\bbeta,\quad i=1,\ldots,n,\:
t=1,\ldots,T.
\end{equation}
The probit version of this model is obtained by assuming that
$\ep_{it}\sim N(0,1)$, so that:
\[
\Phi^{-1}\{p(y_{it}=1|\al_{it},\b x_{it})\} = \al_{it}+\b x_{it}\tr\bbeta,\quad i=1,\ldots,n,\:
t=1,\ldots,T,
\]
with $\Phi(\cdot)^{-1}$
denoting the inverse of the standard Normal cumulative distribution 
function.

Finally, an interesting case (that we consider in the application illustrated
in Section 5) is when each response variable $y_{it}$ is ordinal with categories
$1,\ldots,J$. In this case we can introduce a set of cutpoints $\mu_1\leq\cdots\leq
\mu_{J-1}$ 
and formulate the function in (\ref{eq:link}) as
\[
G(y^*) = \left\{
\begin{array}{cc}
1 & y^*\leq\mu_1,\\
2 & \mu_1<y^*\leq\mu_2,\\
\vdots & \vdots \\
J & y^*>\mu_{J-1}.
\end{array}
\right.
\]
In analogy with the binary case, an ordered logit or an ordered
probit parameterisation results according to whether the error terms
$\ep_{it}$ have a logistic or a standard Normal distribution,
respectively. In the first case, we have that
\begin{equation}\label{eq:global}
\log\frac{p(y_{it}\geq j|\al_{it},\b x_{it})}{p(y_{it}<j|\al_{it},\b x_{it})} =
\mu_j+\al_{it}+\b x_{it}\tr\bbeta,\quad i=1,\ldots,n,\:
t=1,\ldots,T, \: j=2,\ldots,J. 
\end{equation}
This parameterisation is based on global or cumulative logits, the same logits used
in the Odds Proportional model of \cite{mccu:80}.
Note that, as in this model, we are assuming that the effect of the
covariates ($\b x_{it}\tr\bbeta$) and of the unobserved individual
parameters ($\al_{it}$) do not depend on the specific
response category ($j$); this assumption could be removed, but the
model would become more complex to estimate and to interpret.
Finally, note that parameterisation (\ref{eq:global}) is a
generalisation of that in (\ref{eq:logit})
for binary variables.
\subsubsection{Assumption A3}
In order to formulate this assumption, we introduce, for $i=1,\ldots,n$,
the discrete latent variable $u_i$ which has $k$ support points, indexed from 1 to $k$,
and mass probabilities $\pi_1,\ldots,\pi_k$. Then, we assume that
\[
\al_{i1} = \xi_{u_i} + \ep_{i1},\quad i=1,\ldots,n,
\]
and that
\[
\al_{it} = \xi_{u_i} + (\al_{i,t-1} - \xi_{u_i})\rho_{u_i} +
\varepsilon_{it}\sqrt{1-\rho_{u_i}^2},\quad
i=1,\ldots,n,\:t=2,\ldots,T, 
\]
where $\varepsilon_{it}\sim N(0,\sigma^2)$ for all $i$ and $t$. In
the above expressions, $(\xi_h,\rho_h)$ is a pair of parameters
related to the latent state $u_i$ that, for $h=1,\ldots,k$, have to
be estimated jointly with the common variance $\si^2$. In order to
ensure the identifiability of the model, we require that $\xi_1=0$
or, alternatively, $\sum_h \xi_h\pi_h = 0$. The number of parameters
for the latent structure becomes equal to 
$(k-1)+k+(k-1)=3k-2$, which can be
directly compared with those defined for the LAR and LM models in
Section 2.2.

An equivalent way of formulating assumption A3 is by relying, for
$i=1,\ldots,n$, on a standardised AR(1) process of the following
type
\begin{eqnarray*}
\al_{i1}^*&=&\ep^*_{i1},\\
\al_{it}^*&=&\al_{i,t-1}^*\rho_{u_i}+\ep^*_{it}\sqrt{1-\rho_{u_i}^2}, \quad
t=2,\ldots,T,
\end{eqnarray*}
where $\ep^*_{it}\sim N(0,1)$, and then directly including the
parameters $\xi_h$ and $\si^2$ in the equation for the response
variable, that is
\begin{equation}
y_{it}^*=\xi_{u_i}+\al_{it}^*\si+\b x_{it}\tr\b\be+\ep_{it},\quad t=1,\ldots,T.
\label{eq:repara}
\end{equation}
As will be clear in the following, this way to formulate the model is more convenient
for the parameter estimation.
\subsection{Details on latent and manifest distributions}
\begin{figure}[h]\centering
\includegraphics[width=11cm, height=7cm]{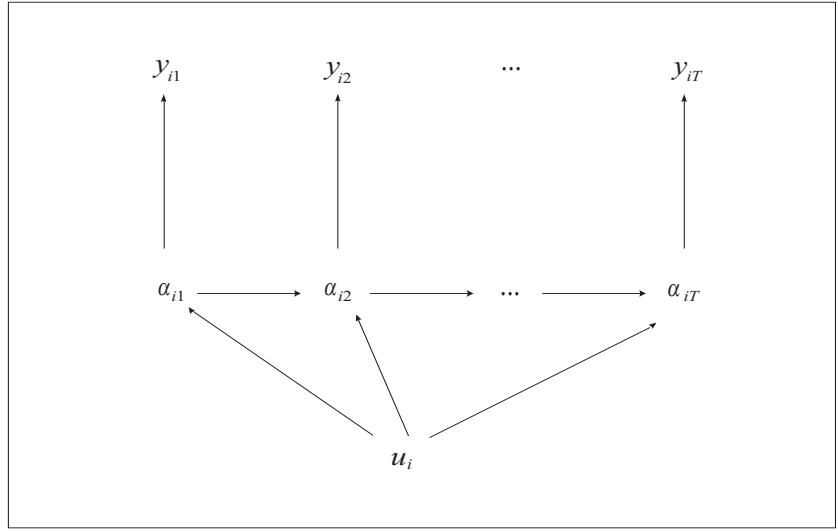}
\caption{\em Path diagram of the 
MLAR model.}
\label{fig:graph1}\vspace*{0.5cm}
\end{figure}

The dependence structure between the latent and
observable variables which results from the above assumptions is
illustrated through the path diagram in Figure \ref{fig:graph1}.

Obviously, the MLAR model generalises the LAR model. In particular,
with $k=1$ the two models coincide, so that in the following they
will be indifferently indicated by LAR or MLAR(1). With $k>1$,
instead, the first model is expected to have a better fit to the
data. This is because the mean value of every $\al_{it}$ and the
correlation coefficient between $\al_{i,t-1}$ and $\al_{it}$ are not
constant, but change 
according to the latent variable $u_i$. We stress that the latent
process based on the above assumptions is still continuous, since
the support of every latent variable $\al_{it}$ is $\Re$. In
particular, we can simply realise that assumption A3 implies the
following mixture model referred to the marginal distribution of
each latent variable:
\[ 
\al_{it}\sim\sum_h\phi(\al_{it};\xi_h,\si^2)\pi_h,\quad i=1,\ldots,n,\:
t=1,\ldots,T,
\]
where $\phi(\al_{it};\xi_h,\si^2)$ is the density function of a
Normal distribution with parameters $\xi_h$ and $\si^2$. Similarly,
concerning the marginal distribution of $(\al_{i,t-1},\al_{it})$ we
have that
\[
(\al_{i,t-1}, \al_{it})\tr\sim\sum_h\phi_2((\al_{i, t-1},\al_{it})\tr; \xi_h\b 1_2, \bSigma_{u_i})\pi_h,
\]
which now involves the density function of a bivariate Normal
distribution with mean $\xi_h\b 1_2$, where $\b 1_2$ denote a vector of two
ones,
and variance-covariance matrix
\begin{displaymath}
\bSigma_h = \sigma^2
\left( \begin{array}{cc}
1 & \rho_h \\
\rho_h & 1  \\
\end{array} \right).
\end{displaymath}

The above arguments 
imply that a possible interpretation of the MLAR
model may be based on considering the population of subjects, from
which the observed sample comes, as made of $k$ subpopulations (or
latent classes), such that a LAR model with the same parameters
holds within each subpopulation. In fact, the probability mass (or
density) function of the distribution of the response vector $\b
y_i$ given all the observable covariates $\b X_i$ may be expressed
as a {\em mixture of LAR models}. In particular, we have the
following manifest distribution:
\begin{equation}
p(\b y_i|\b X_i)=\sum_h p^{(h)}(\b y_i|\b X_i)\pi_h,\label{eq:manifest}
\end{equation}
with $p^{(h)}(\b y_i|\b X_i)$ defined as in (\ref{eq:mani_lar}), 
for $h = 1, \ldots, k$.
However, in order to implement the estimation method for the model
parameters, it is more convenient to express this probability 
or density 
on the basis of the latent effect
$\al_{it}^*$ which follows a standardised AR(1) process. Then, we
have
\begin{eqnarray}
p^{(h)}(\b y_i|\b X_i)&=& \int p(y_{i1}|\al_{i1}^*,\b
x_{i1})f^{(h)}(\al_{i1}^*)
\int p(y_{i2}|\al_{i2}^*,\b x_{i2})f^{(h)}(\al_{i2}^*|\al_{i1}^*)\cdots\nonumber\\
&&\cdots\int p(y_{iT}|\al_{iT}^*,\b
x_{iT})f^{(h)}(\al_{iT}^*|\al_{i,T-1}^*) d\al_{iT}\cdots
d\al_{i2}d\al_{i1},\label{eq:pcond}
\end{eqnarray}
where $p(y_{it}|\al_{it}^*,\b x_{it})$ is computed on the basis of
(\ref{eq:repara}) and
\begin{eqnarray}
f^{(h)}(\al_{i1}^*)&=&\phi(\al_{i1}^*;0,1),\label{eq:dens1}\\
f^{(h)}(\al_{it}^*|\al_{i,t-1}^*)&=&\phi(\al_{it}^*;\al_{i,t-1}^*\rho_h,1-\rho^2_h).
\label{eq:dens2}
\end{eqnarray}
An efficient way to compute 
function (\ref{eq:pcond}) is described in Section 4.1.

However, our main aim here is not that of classifying subjects in
different subpopulations, but that of having a flexible structure for the latent
process. In fact, by rising $k$, we have an increasing 
degree of flexibility of the distribution of $\b\al_i$ with respect
to assuming a standard AR(1) process as in the LAR model. In fact,
it is well-known that, with a suitable number of components and under suitable
conditions, a mixture distribution can adequately approximate any distribution.
The same principle has been exploited by \cite{Bart:clus:2005} to
propose a flexible method to classify univariate observations and by
\cite{Scac:Bart:hier:2005} to propose a regression model with a
flexible distribution for the error terms.

\begin{figure}[!h]\centering
\begin{tabular}{cc}
\includegraphics[width=8cm]{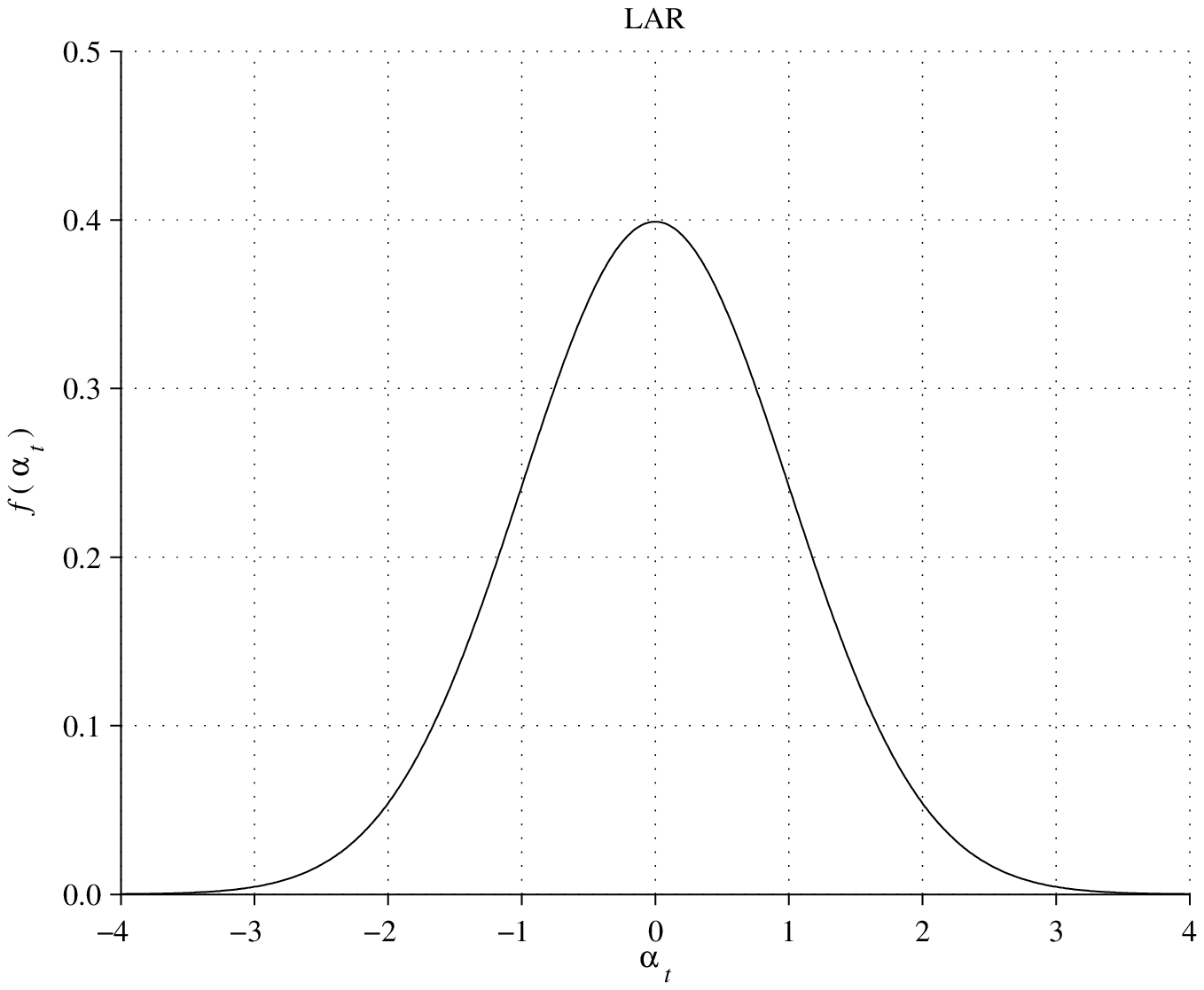} &
\includegraphics[width=8cm]{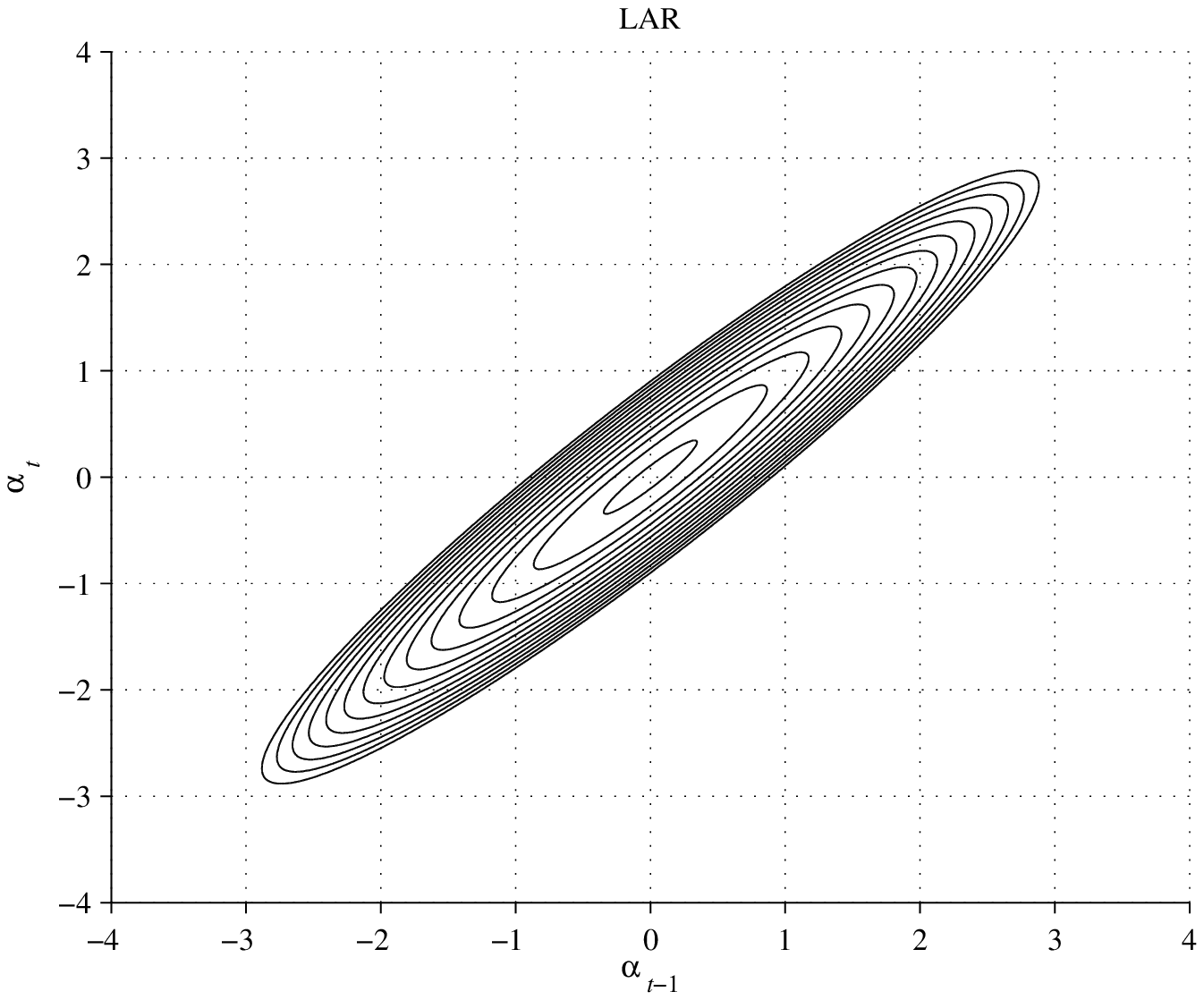}\\
\includegraphics[width=8cm]{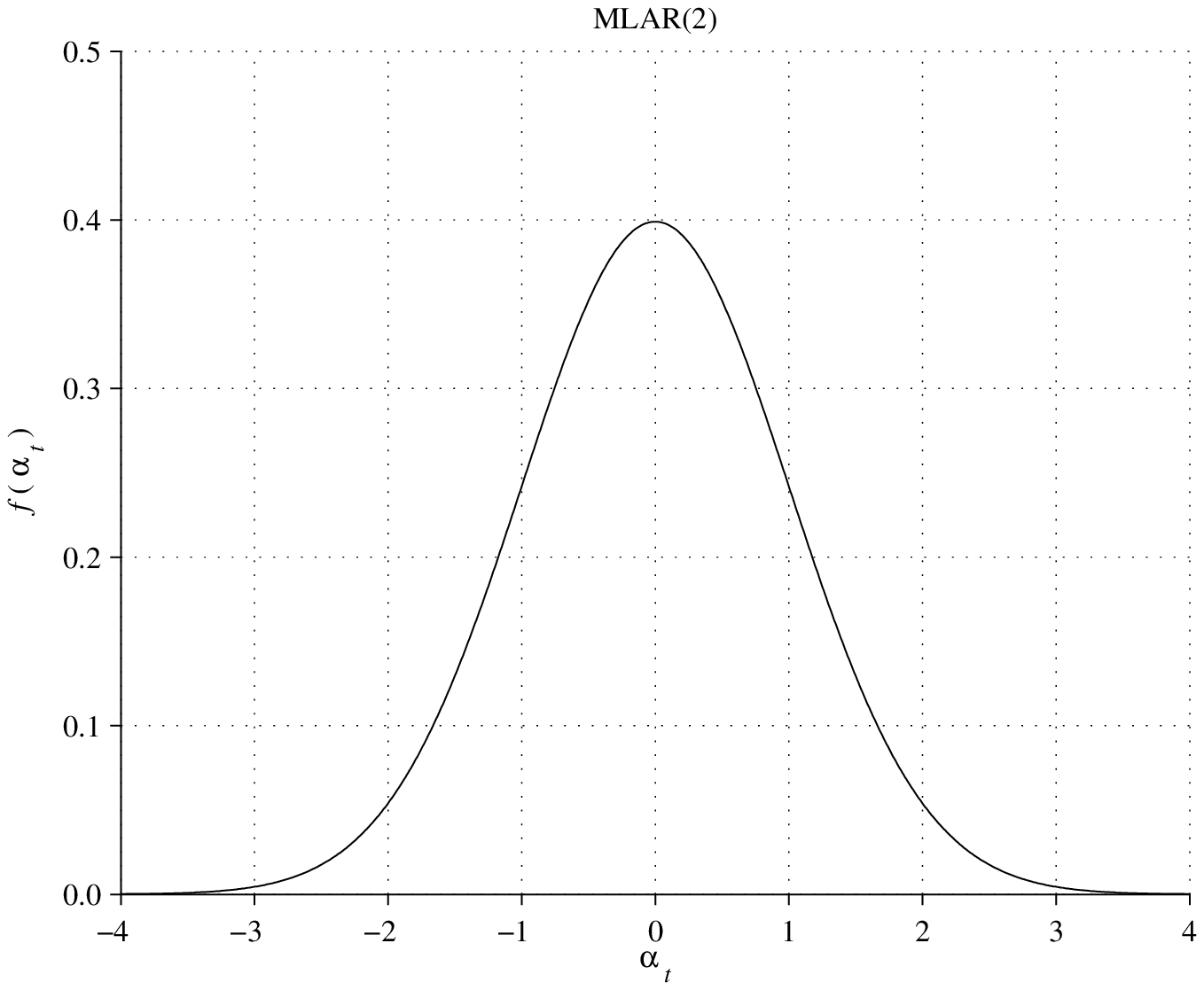} &
\includegraphics[width=8cm]{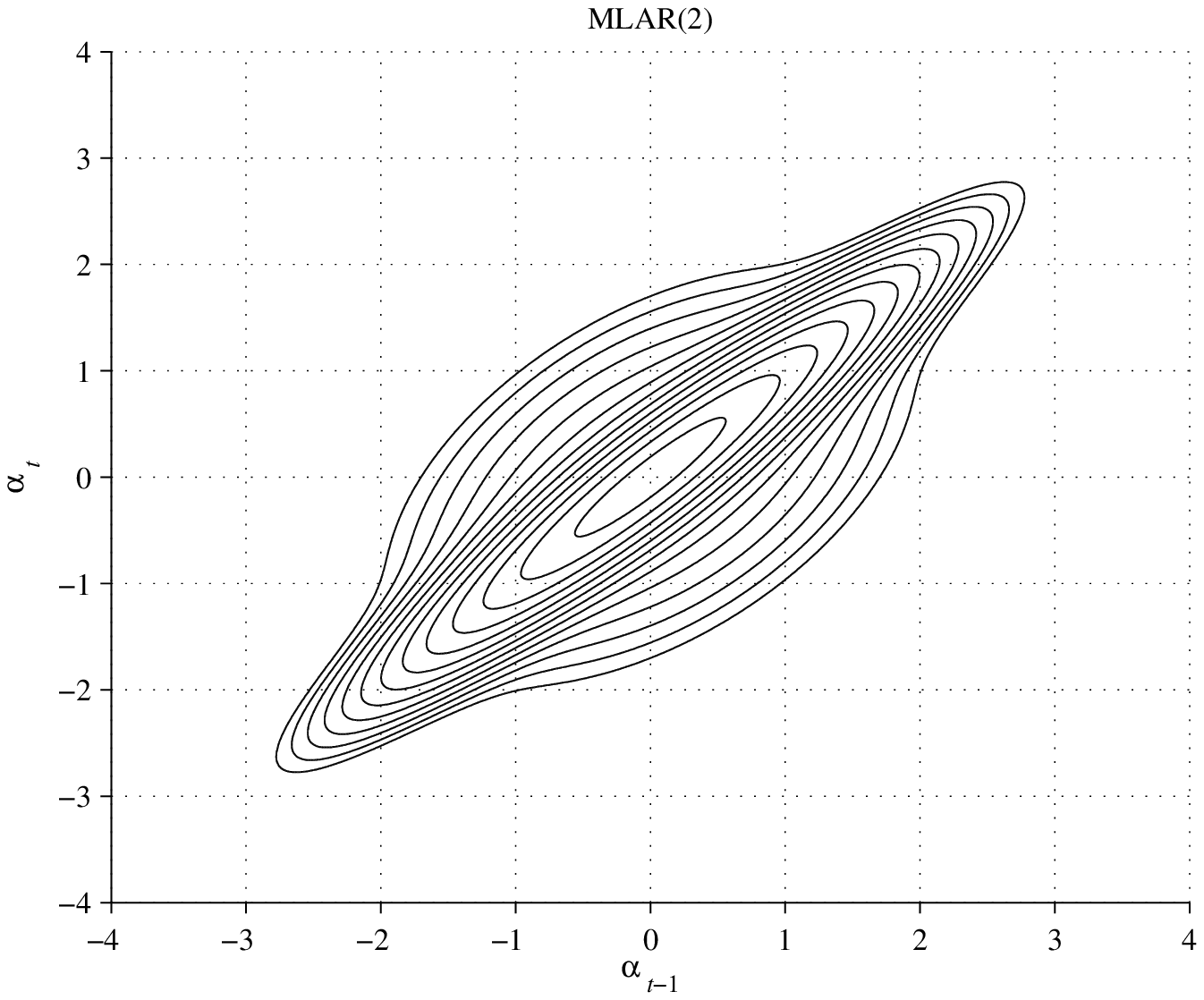}\\
\includegraphics[width=8cm]{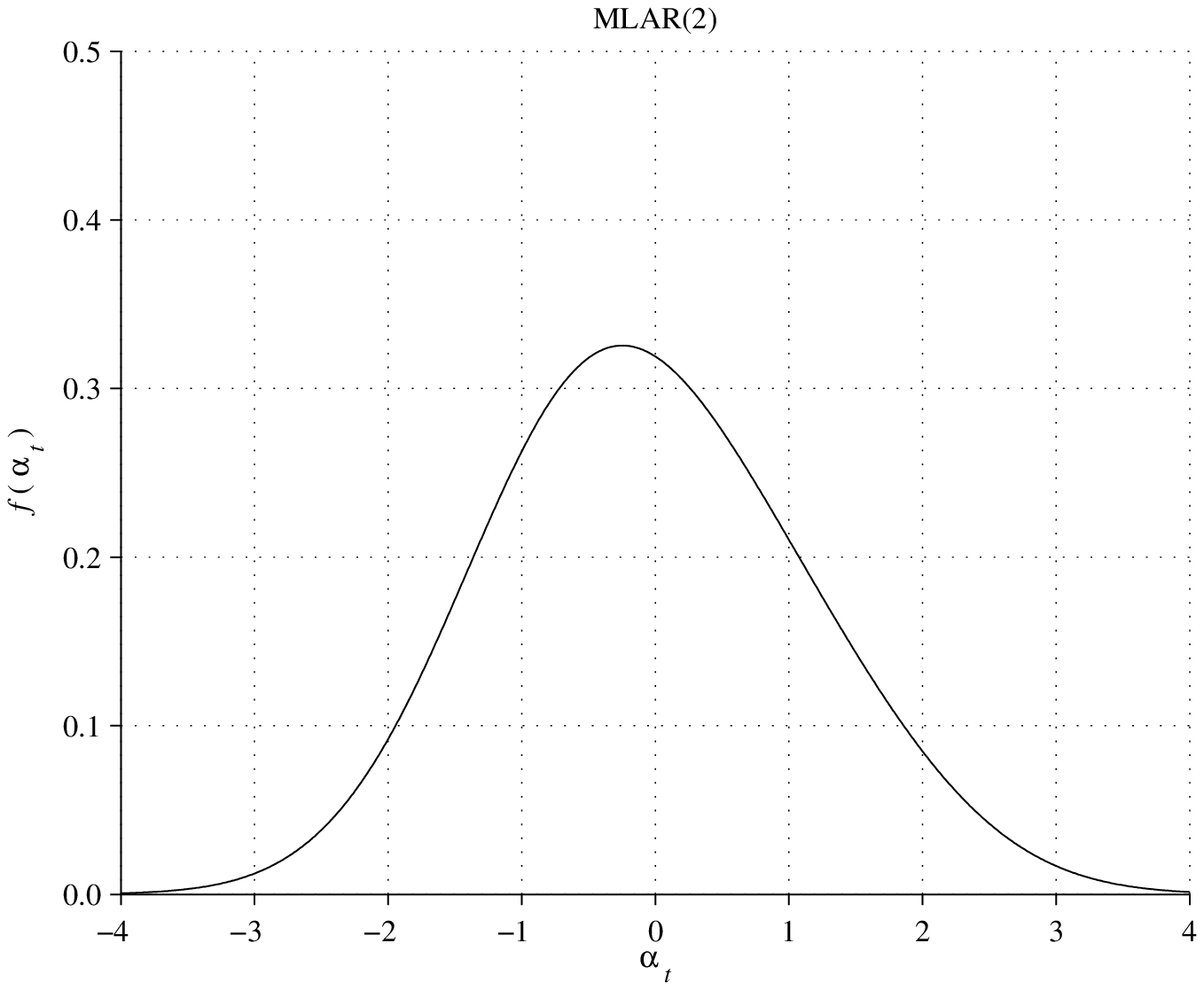} &
\includegraphics[width=8cm]{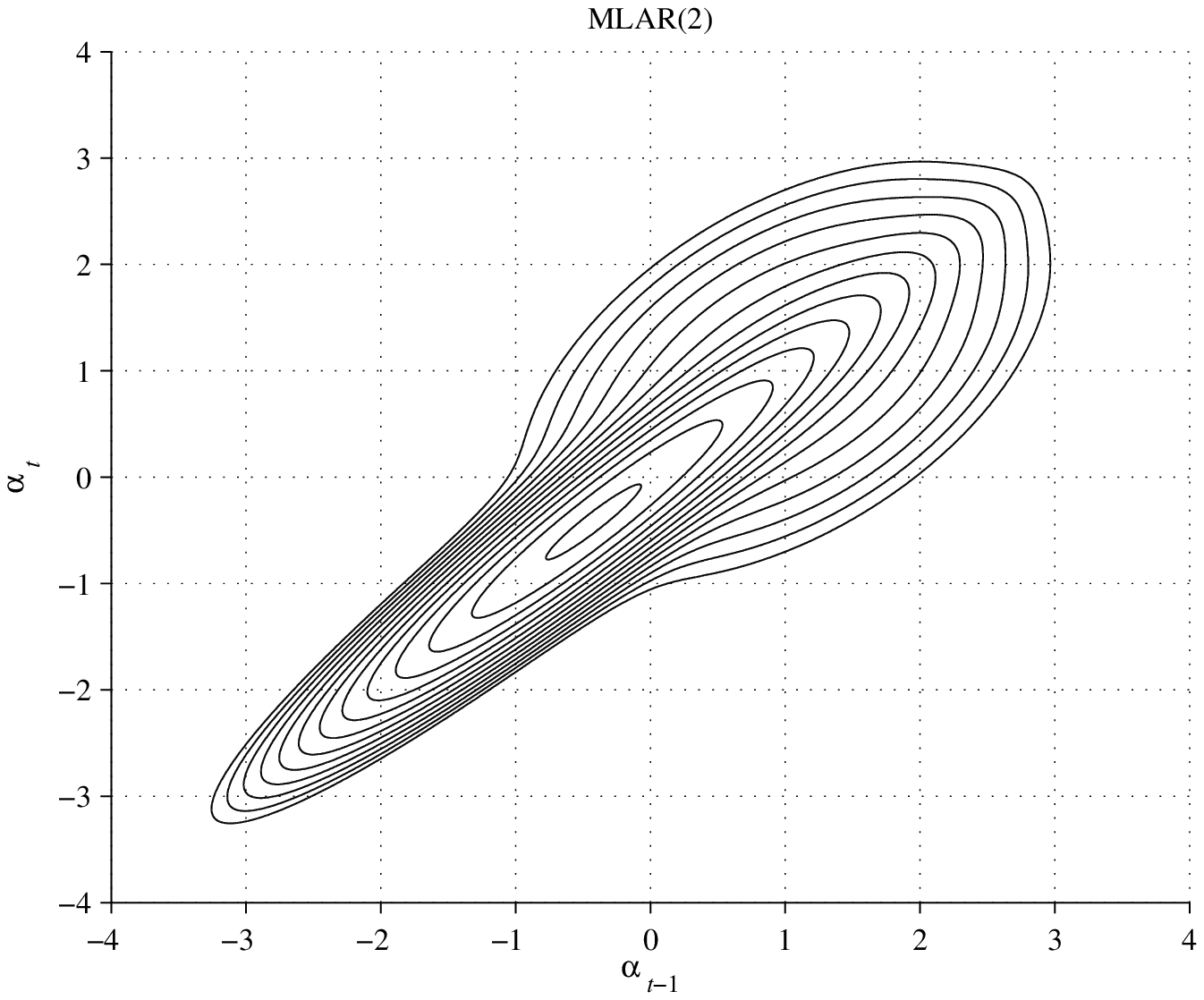}\\
\end{tabular}
\caption{\em Density function of the univariate distribution of $\al_{it}$
(left) and of the bivariate distribution of $(\al_{i,t-1},\al_{it})$
(right) under the LAR model (top) with 
$\rho=0.95$ and $\si^2=1.00$
and under the MLAR(2) model with 
$\rho_1=0.95$, $\rho_2=0.50$,
$\pi_1=0.70$, $\pi_2=0.30$, and $\si^2=1.00$
and with $\xi_1=\xi_2=0.00$ (middle)
and $\xi_1=-0.50$ and $\xi_2=1.00$ (bottom).}
\label{fig:teo}
\end{figure}

In order to clarify the above point, in Figure \ref{fig:teo} we
represent the density function of the marginal distribution of every
$\al_{it}$ and of $(\al_{i, t-1},\al_{it})$ for the LAR model and
for a MLAR model with $k=2$ components and different parameter
values. In particular, the top panel in Figure \ref{fig:teo} is
referred to LAR model with parameters $\rho=0.95$ and $\si^2=1.00$,
the middle panel is referred to the MLAR(2) model with parameters
$\xi_1=\xi_2=0.00$, $\rho_1=0.95$, $\rho_2=0.50$, $\pi_1=0.70$ and
$\pi_2=0.30$, whereas the bottom panel is referred to the same
MLAR(2) model with $\xi_1=-0.50$ and $\xi_2=1.00$. In order to make
clear the comparison among the plots, 
we used the same level curves for
all bivariate distributions in Figure \ref{fig:teo} (right panels),
which are defined on the basis of a grid of equispaced points 
on the logarithmic scale.

We observe that with only two components, very different shapes of 
the density function 
of the latent variable distribution 
may be obtained. In particular, when $\xi_1=\xi_2$ (middle panel of Figure
\ref{fig:teo}), both univariate and bivariate distributions are
still symmetric. However, the bivariate density function has a
different shape with respect to that under the LAR model, due to a
much higher dispersion around the middle of the plot. Moreover, with
$\xi_1\neq\xi_2$ (bottom panel), these distributions are also
asymmetric, with the density of points in the North-West region of
the bivariate plot that considerably rises. In a similar way we can
even generate more complex shapes if we use, for instance, values of
$k$ higher than 2, at the cost of a moderate increase of the number
of parameters.
\section{Likelihood inference}
In this section, we deal with likelihood inference for the model
proposed in Section 3. In particular, we first show how to
efficiently compute the model log-likelihood. Then we deal with its
maximisation by an EM algorithm and we describe how to compute
standard errors and predict individual effects. Finally, we deal
with the choice of the number of mixture components.
\subsection{Computation of the model likelihood}
Since the sample units are assumed to be independent, the model likelihood
has logarithm
\begin{equation}
\ell(\b\th) = \sum_i\log p(\b y_i|\b X_i),\label{eq:log-lik}
\end{equation}
where $\b\th$ is a short hand notation for all the non-redundant
model parameters and, for every subject $i$, $p(\b y_i|\b X_i)$
denotes the probability mass (or density) function given in
(\ref{eq:manifest}), seen as a function of these parameters.

In order to efficiently compute $p(\b y_i|\b X_i)$, we exploit a
recursion developed in the hidden Markov literature. First of all,
we transform the series of integrals to compute $p^{(h)}(\b y_i|\b
X_i)$, which is defined in (\ref{eq:pcond}), in a series of sums on
a suitable grid of quadrature points, as
proposed by \cite{hei:08}. Let $q$ denote the number of quadrature
points and let $\nu_m$ denote the $m$-th quadrature knot, with
$m=1,\ldots,q$. 
Moreover, for each mixture component $h$, 
let $\omega_m^{(h)}$ denote the $m$-th weight for
the integral with respect to $\al_{i1}$ and let
$\omega_{m_1m_2}^{(h)}$ denote the $m_2$-th weight for the integral
with respect to $\al_{it}$, given
that the $m_1$-th knot is selected for the integral with respect to
$\al_{i,t-1}$. 
Then, the expression in (\ref{eq:pcond}) becomes:
\begin{equation}
p^{(h)}(\b y_i|\b X_i)=\sum_{m_1} p(y_{i1}|\nu_{m_1},\b
x_{i1})\omega^{(h)}_{m_1} \sum_{m_2} p(y_{i2}|\nu_{m_2},\b
x_{i2})\omega^{(h)}_{m_1m_2}\cdots \sum_{m_T}p(y_{iT}|\nu_{m_T},\b
x_{iT})\omega^{(h)}_{m_{T-1}m_T}.\label{eq:quadrature}
\end{equation}
In practice, the knots are taken on a suitable grid of points
between two extremes, say $-5$ and $5$, and the corresponding
weights are computed as follows:
\begin{eqnarray}
\omega_m^{(h)}&=&\frac{f^{(h)}(\nu_m)}
{\sum_lf^{(h)}(\nu_l)},\quad m=1,\ldots,q,\nonumber\\
\omega_{m_1m_2}^{(h)}&=&
\frac{f^{(h)}(\nu_{m_2}|\nu_{m_1})}
{\sum_lf^{(h)}(\nu_{m_2}|\nu_l)},\quad m_1,m_2=1,\ldots,q,\label{eq:bivweight}
\end{eqnarray}
for $h=1,\ldots,k$, on the basis of the density functions defined in
(\ref{eq:dens1}) and (\ref{eq:dens2}). The quadrature knots and the
corresponding weights could be found by a more complex method, such
as the Guass-Hermite method. However, we experienced that the above
method leads to essentially equivalent solutions when $q$ is large
enough.

We can easily recognise that expression (\ref{eq:quadrature}) is the
same expression of the manifest distribution of the LM model given
in (\ref{eq:manifest_LM}). The only difference is that for the LM
model we have support points ($\xi_h$) and initial and transition
probabilities ($\pi_h$, $\pi_{h_1h_2}$) to be estimated on the basis
of the data. Here, we have knots ($\nu_m$) and weights
($\omega^{(h)}_m$, $\omega^{(h)}_{m_1m_2}$) which are instead given,
with the exception of the weights $\omega^{(h)}_{m_1m_2}$ which only
depend on the correlation coefficient $\rho_h$. However, the same
recursion of \cite{bau:70} may be used to efficiently compute
(\ref{eq:quadrature}), and then obtain $p^{(h)}(\b y_i|\b X_i)$ from
which we obtain $p(\b y_i|\b X_i)$ by (\ref{eq:manifest}) and the
log-likelihood $\ell(\b\th)$ by (\ref{eq:log-lik}). Note that
applying the recursion at issue is essentially equivalent to apply
the SGQ of \cite{hei:08}
and to the method of \cite{bart:delu:01,bart:delu:03} to compute the
likelihood function of SV models.
\subsection{Maximum likelihood estimation}
As derived above, once a suitable set of quadrature knots has been
adopted, the likelihood of the proposed model may be seen as
equivalent to that of an LM model
with covariates and latent parameters suitably
constrained. Then, maximum likelihood (ML) estimation may be
performed by an adaptation of the EM algorithm for the LM model
described by \cite{bart:farc:09}; see also \cite{bau:70} and
\cite{demp:lair:rubi:77}. In the following, we outline this extended
algorithm, referring for some details to \cite{bart:farc:09}.

The EM algorithm is based on the so-called {\em complete data log-likelihood}
that, in the present case, corresponds to the log-likelihood that we could compute
if we knew, for $i=1,\ldots,n$, the value of the latent variable $u_i$ and the
value of the quadrature knot for $\al_{it}$, $t=1,\ldots,T$. This is equivalent to
the knowledge of the dummy variables $w_{ih}$, $i=1,\ldots,n$, $h=1,\ldots,k$,
and $z_{imt}$, $i=1,\ldots,n$,
$m=1,\ldots,q$, $t=1,\ldots,T$, where $w_{ih}=I\{u_i=h\}$ and $z_{imt}=
I\{\al_{it}=\nu_m\}$. Up to a constant term,
the complete data log-likelihood may be expressed as
\begin{equation}
\ell^*(\b\th) = \sum_i\sum_h w_{ih}\left\{\log\pi_h+
\sum_{m_1}\sum_{m_2}\sum_{t>1} z^*_{im_1m_2t}\log\omega^{(h)}_{m_1m_2}+
\sum_m\sum_t z_{imt}\log p(y_{it}|\nu_m,\b x_{it})\right\},\label{eq:comp_lk}
\end{equation}
where $z^*_{im_1m_2t}=z_{im_1,t-1}z_{im_2t}$.

The EM algorithm alternates the following steps until convergence:
\begin{itemize}
\item {\bf E-step}: compute the conditional expected value of the complete data
log-likelihood given the observed data and the
current estimate of $\b\th$;
\item {\bf M-step}: maximise the expected value above with respect to $\b\th$.
\end{itemize}

The E-step is equivalent to computing the conditional expected value, given the
observed data, of every dummy variable $w_{ih}$ and of the products $w_{ih}z_{imt}$
and $w_{ih}z_{im_1m_2t}^*$.
In practice, for $i=1,\ldots,n$, we have that
\[
\hat{w}_{ih}=E(w_{ih}|data) = \frac{p^{(h)}(\b y_i|\b X_i)}{p(\b y_i|\b X_i)},\quad
h=1,\ldots,k,
\]
where the probabilities are computed on the basis of the current
value of the parameters and {\em data} stands for ``observed data''.
Moreover, we have that
\[
\widehat{(w_{ih}z_{imt})}=
E(w_{ih}z_{imt}|data) = \hat{w}_{ih}E(z_{imt}|w_{ih}=1,data),\quad m=1,\ldots,q,
\: t=1,\ldots,T,
\]
where $E(z_{imt}|w_{ih}=1,data)$ is the {\em posterior probability} that
subject $i$ is in state $h$ at time occasion $t$ given that $u_i=h$, and
\[
\widehat{(w_{ih}z^*_{im_1m_2t})}=
E(w_{ih}z^*_{im_1m_2t}|data) = \hat{w}_{ih}E(z^*_{im_1m_2t}|w_{ih}=1,data),
\quad m_1,m_2=1,\ldots,q,\: t=2,\ldots,T,
\]
where $E(z^*_{im_1m_2t}|w_{ih}=1,data)$ is the posterior probability that
subject $i$ moves from state $m_1$ to state $m_2$ at occasion $t$, given that $u_i=h$.
These posterior probabilities may be computed by suitable recursions; see
\cite{bau:70}, \cite{bart:farc:09}, and \cite{bart:farc:penn:10} for details.

Once the expected values of the dummy variables have been
substituted in (\ref{eq:comp_lk}), the resulting function is
maximised with respect to the model parameters, which are
consequently updated. The easiest parameters to update are the
probabilities $\pi_h$, for which we have an explicit solution
\[
\pi_h = \frac{\sum_i\hat{w}_{hi}}
{\sum_l\sum_i\hat{w}_{li}},\quad h=1,\ldots,k.
\]
Then,
in order to update each parameter $\rho_h$, $h=1,\ldots,k$, we
have to maximise, by a numerical optimisation algorithm, the
function
\[
\sum_i\sum_{m_1}\sum_{m_2}\sum_{t>1}
\widehat{(w_{ih}z^*_{im_1m_2t})}\log\omega^{(h)}_{m_1m_2},
\]
which depends on this parameter through (\ref{eq:bivweight}).
Finally, the other model parameters, that is $\xi_1,\ldots,\xi_k$,
$\b\be$, and $\sigma^2$, are update by maximising the function
\[
\sum_i\sum_h
\sum_m\sum_t \widehat{(w_{ih}z_{imt})}\log p(y_{it}|\nu_m,\b x_{it}),
\]
which depends on these parameters through (\ref{eq:repara}). This maximisation may
be performed by a NR iterative algorithm, the implementation of which
is not difficult, due to the availability of explicit expressions for the
first and second derivatives of the target function.

Since the EM algorithm is rather slow to converge, after a certain number of steps
we switch to a full NR algorithm to maximise the model log-likelihood
$\ell(\b\th)$. This algorithm updates the model parameters $\b\th$ by adding the
following quantity $\b J(\b\th)^{-1}\b s(\b\th)$, where $\b s(\b\th)$ denotes the
score vector for $\ell(\b\th)$ and $\b J(\b\th)$ denotes the corresponding
observed information
matrix. The latter is equal to minus the second derivative of $\ell(\b\th)$ with respect
to $\b\th$. Following \cite{bart:farc:09}, the score vector is computed as the
first derivative of the expected value of complete data log-likelihood, which is
obtained after an E-step. The observed information matrix is then obtained on the basis
of the numerical derivative of $\b s(\b\th)$.

We take the value of $\b\th$ at convergence of the NR algorithm as
the ML estimate $\hat{\b\th}$. As it typically happens for
latent variable models, the model likelihood may be multimodal and the
point at convergence depends on the starting values for the
parameters, which need to be carefully chosen. Then, we suggest to try different
starting values in order to be sure that the found solution corresponds to the
global maximum of $\ell(\b\th)$.

Once the ML estimates have been computed, it may be of interest to
obtain the corresponding standard errors. These may be obtained in
the usual way on the basis of $\b J(\hat{\b\th})^{-1}$ and may be
used to compute confidence intervals for the parameters and perform
Wald testing about certain hypotheses of interest. More generally,
hypotheses of interest may be tested by a likelihood ratio statistic
that, under the usual regularity conditions, has asymptotic
distribution of $\chi^2$-type.

On the basis of the parameter estimates it may also be of interest
to predict every latent variable $\al_{it}$. This may be performed
through the following posterior expected value given the observed
data:
\begin{equation}
\hat{\al}_{it}=\sum_h\sum_m\widehat{(w_{ih}z_{imt})}(\hat{\xi}_h+\nu_m\hat{\sigma}),
\quad i=1,\ldots,n,\:t=1,\ldots,T,
\label{eq:pred}
\end{equation}
with all quantities computed on the basis of the final estimate $\hat{\b\th}$.

For the case of binary and ordinal response variables, we
implemented the above strategy to obtain the ML estimate of $\b\th$,
which is based on the joint use of an EM and of an NR algorithm, in
a series of {\sc Matlab} functions that we make available to the
reader upon request. In our experience, this strategy 
properly works
and provides ML estimates and corresponding standard errors in a
reasonable amount of time, provided that $k$ is not too large. 
%
\subsection{Selection of the number of mixture components and of
quadrature points}\label{sec:model_selection}
In applying the proposed model, of crucial importance is the choice
of the number of mixture components ($k$) and of quadrature points
($q$). Regarding the choice of $q$, we have to use a value which is
large enough to guarantee an adequate approximation of the true
likelihood function, that is the likelihood that we could obtain by
exactly computing the multiple integral in (\ref{eq:pcond}). At this
regard, the strategy we suggest is based on trying, for a given $k$,
increasing values of $q$ until the maximum of $\ell(\b\th)$ does not
significantly change with respect to the previous value of $q$. In
our application, for instance, we start 
with $q=21$ and we increase 
the value of $q$ by 10 at each attempt, stopping when the
maximum of $\ell(\b\th)$ increases less than $10^{-3}$.

Through the above strategy, we find 
a suitable value of
$q$ for a given $k$. The point now is how to choose $k$. In summary,
the strategy we suggest consists of increasing $k$ until the
estimated latent structure does not significantly change. In
practice, for each tried value of $k$ we obtain the predicted latent
variables $\al_{it}$ through (\ref{eq:pred}) on the basis of
$\hat{\b\th}$ and, for $k>1$, we compute the correlation index
between these predicted values and those computed with $k-1$ mixture
components. The first value of $k$ such that this correlation index
is higher than a suitable threshold (we use 0.99 in our application)
is taken as the optimal number of mixture components.

Note that, in order to select $k$, we could also rely on information
criterion such as the Akaike Information Criterion \citep{aka:73} or
the Bayesian Information Criterion \citep{sch:78}. However, we
experimented in our applications that these criteria tend to choose
a value of $k$ higher than necessary, whereas we have evidence that
the criterion suggested above, which is based on direct assessment
of the estimated latent structure, has good performance.
\section{Application to Self-reported health status}
To illustrate the proposed approach, we consider a dataset which derives from
the Health and Retirement Study conducted by the University of Michigan (see
\verb"http://www.rand.org/labor/aging/dataprod" for detailed illustration).
After a description of the dataset, we report the results of its analysis based
on the proposed approach.
\subsection{Dataset description}
The dataset is referred to a sample of $n=7,074$ American individuals
who were 
asked to express opinions on their health
status at $T=8$ approximately equally spaced occasions, from 1992 to 2006.
The response variable ({\em self-reported health status}) is measured on a scale
based on five ordered categories:
``poor'', ``fair'', ``good'', ``very good'', and
``excellent''. For every subject some covariates are also
available: {\em gender}, {\em race}, {\em education},
and {\em age} (at each time
occasion). Table \ref{tab:descr} shows some descriptive statistics about these
covariates, whereas Table \ref{tab:descr1} shows the marginal
distribution of the response variable over the 8
occasions of interview.

\begin{table}[h!]\centering
\vspace*{0.25cm}
{\small
\begin{tabular}{ll|ccc}
\hline \hline
Variable    &   Category &  $\%$  & Mean & St.Dev.    \\
\hline
\emph{gender}:   &   female  &   58.1    &   --   &  -- \\
    &    male       &    41.9   &   --    &    --  \\
\hline
\emph{race}:       &  white       &      82.9    &  --   &  --   \\
    &   non white   &   17.1    &   --  & --    \\
\hline
\emph{education}:      &  high school  &   60.9    &   --    & -- \\
    &  some college     &   19.7    &   --  & -- \\
    &  college and above    &   19.4    &   --  &   -- \\
\hline
 \emph{age} (in 1992): &    &     -- &  54.8     &  5.5   \\
\hline \hline
\end{tabular}}
\caption{\em Distribution of the covariates. \label{tab:descr}}
\vspace*{0.25cm}\end{table}

\begin{table}[h!]\centering
\vspace*{0.25cm}
{\small
\begin{tabular}{l|rrrrrrrr|r}
\hline \hline
    & \multicolumn{8}{c|}{occasion of interview}  &  \\
\hline
SRH category    &   \multicolumn1c{1}   &   \multicolumn1c{2}   &
\multicolumn1c{3}   &   \multicolumn1c{4}   &   \multicolumn1c{5}   &
\multicolumn1c{6}   &   \multicolumn1c{7}   & \multicolumn1{c|}{8}  &  Total\\
\hline
poor    &   4.7 &   4.7 &   4.7 &   6.1 &   5.5 &   5.9 &   7.2 &   8.1 &   5.9 \\
fair    &   11.5    &   13.0    &   13.1    &   17.0    &   15.6    &   16.9    &
19.2    &   20.3    &   15.8    \\
good    &   27.2    &   28.9    &   28.1    &   32.0    &   30.6    &   32.0    &
32.2    &   32.0    &   30.4    \\
very good   &   30.8    &   32.6    &   34.3    &   31.1    &   33.6    &   32.3    &
30.0    &   29.7    &   31.8    \\
excellent   &   25.7    &   20.8    &   19.8    &   13.7    &   14.7    &   13.0    &
11.4    &   10.0    &   16.1    \\
\hline
Total   &   100.0   &   100.0   &   100.0   &   100.0   &   100.0   &   100.0   &
100.0   &   100.0   &   100.0   \\
\hline \hline
\end{tabular}}
\caption{\em Distribution of the response variable over the occasions of interview
in terms of percentage frequencies. \label{tab:descr1}}
\vspace*{0.25cm}
\end{table}

As shown in Table \ref{tab:descr}, the main part of individuals
in the sample are females ($58.1\%$)
and whites ($82.9\%$), with an average age at the first time occasion equal to $54.8$;
we recall that the occasions of interview are around two years far apart.
The $60.9\%$ of the sample has a high-school
diploma, whereas a college degree or a higher title is possessed by the
$19.4\%$ of subjects. In the following, the covariate \emph{education} is
introduced in the model by
assigning increasing scores to its categories: 1 for ``high school'', 2 for
``some college'' (i.e., a high school or a general education diploma and more than 12
years of education), 3 for ``college and above''
(i.e., a college degree, such as Bachelor of Arts,
or an higher title, such as PhD).

About the distribution of the response variable at each time occasion
(see Table \ref{tab:descr1}),
we observe that more than the $60\%$ of responses is equally distributed between
categories  ``good'' and ``very good'', being substantially stable over time.
Moreover, the $16.1\%$ of
individuals evaluates the health status as ``excellent'', with a decreasing trend
(the percentage is over $25.7\%$ at the first occasion and it decreases to
$10.0\%$ at the eighth one).
On the other side, the remaining part of individuals gives a
negative  judgement to the health status, with increasing percentages over time: from
$4.7\%$ to $8.1\%$ and from $11.5\%$ to $20.3\%$ for categories
``poor'' and ``fair'' response, respectively.

More insights about the 
subjects' responses to the questionnaire may be
derived on the basis of the empirical transition matrix reported in
Table \ref{tab:descr3}. Each row of this matrix shows the percentage
frequencies of the five response categories at occasion $t$ given the response
at occasion $t-1$, with $t=2,\ldots,T$.

\begin{table}[!ht]\centering
\vspace*{0.25cm}
{\small
\begin{tabular}{l|rrrrr|r}
\hline \hline
   & \multicolumn5c{SRH at $t$}   \\
\hline
SRH at $t-1$    &   \;\;\;\;\;poor      & \;\;\;\;\;fair        &   \;\;\;\;\;good  &
very good   &   excellent   & Total \\
\hline
poor    &   54.5    &   34.1    &   8.4 &   2.5 &   0.7  &  100.0   \\
fair    &   12.8    &   51.0    &   27.4    &   7.2 &   1.6  &  100.0\\
good    &   2.5 &   16.5    &   53.3    &   23.6    &   4.1 &  100.0    \\
very good   &   0.8 &   4.7 &   25.9    &   55.6    &   13.0    &  100.0 \\
excellent   &   0.4 &   1.9 &   10.6    &   33.7    &   53.4 &  100.0   \\
\hline \hline
\end{tabular}}
\caption{\em Conditional empirical distributions of the response variable
at time $t$ given the response at time $t-1$,
with $t=2,\ldots,T$ (percentage frequencies). \label{tab:descr3}}
\vspace*{0.25cm}
\end{table}

In general, a rather high persistence of the judgement about the health status results,
since more than one half of responses at time
$t$ is in the same category as 
the response at time $t-1$, and percentages included
between $12.8\%$ and $34.0\%$ lie in an adjacent category. On the other hand, jumps
between different and not adjacent response categories in consecutive time occasions are
observable for the remaining part of the sample. For example, among subjects who
evaluate their health status as ``poor'' at a given occasion,
the $8.4\%$ evaluates 
it as ``good'' at the next
occasion; on the contrary, among subjects who respond 
``very good''
at a given occasion, only the $4.7\%$ responds ``fair'' 
at the following occasion.
\subsection{Model selection}\label{sec:application:model_selection}
To the data described above, we preliminary fit the proposed MLAR model for
different values of $k$ (number of mixture components) and $q$ (number of quadrature
points). To take into account the ordinal nature of the response variable, the
model is formulated on the basis of the global logit parameterisation defined
in equation (\ref{eq:global}). Then, the optimal values of $k$ and $q$ is chosen as
described in Section \ref{sec:model_selection}. We recall that, as concerns the
selection of $q$ (given $k$)
the adopted procedure starts from $q = 21$ and increases it by 10;
the number
of quadrature points is selected in correspondence of the first difference
between two consecutive maximum log-likelihood 
values smaller than 0.001.
With reference to the selection of $k$, the adopted strategy consists of computing
the correlation index $\rho_{k-1,k}$ between the predicted
$\al_{it}$ values of MLAR($k$) and those of MLAR($k-1$), for
increasing values of $k$ starting from 2.
When $\rho_{k-1,k}$ is greater than 0.99, $k$ is not raised anymore and
its last value is taken as the optimal number of mixture components.
Table \ref{Table_choice_kq} reports the main results of this model
selection procedure.

\begin{table}[!ht]\centering
{\small
\begin{tabular}{r|rr|rr|rr}
\hline \hline
    &   \multicolumn{2}{c|}{$k=1$}          &   \multicolumn{2}{c|}{$k=2$}      &   \multicolumn{2}{c}{$k=3$}       \\
\multicolumn{1}{c|}{$q$} &   log-likelihood  &   difference  &   log-likelihood  &
 difference  &   log-likelihood  &   difference  \\
\hline
21  &   -63609.195  &   --  &   -62978.009  &  --   &   -62831.584  &   --  \\
31  &   -63624.648  &   -15.453 &   -62996.639  &   -18.630 &   -62844.763  &   -13.179 \\
41  &   -63624.657  &   -0.009  &   -62998.591  &   -1.952  &   -62845.683  &   -0.920  \\
51  &   {\bf -63624.657}    &   \textbf{0.000}  &   -62998.613  &   -0.022  &   -62845.688  &   -0.005  \\
61  &   -63624.657  &   0.000   &   {\bf -62998.613}    &   \textbf{0.000}
    &   {\bf -62845.688}    &   \textbf{0.000}  \\
71  &   -63624.657  &   0.000   &   -62998.613  &   0.000   &   -62845.688  &   0.000   \\
81  &   -63624.658  &   0.000   &   -62998.614  &   0.000   &   -62845.688  &   0.000   \\
91  &   -63624.658  &   0.000   &   -62998.614  &   0.000   &   -62845.688  &   0.000   \\
101 &   -63624.658  &   0.000   &   -62998.614  &   0.000   &   -62845.688  &   0.000   \\
\hline
&      &           \multicolumn2c{$\rho_{12}=0.9783$}  &
\multicolumn2c{$\rho_{23}=0.9974$} &            \\
\hline \hline
\end{tabular}}
\caption{\em Log-likelihoods and differences between consecutive values for $k=1,2,3$
and $q$ from $21$ to $101$ with 
step 10; $\rho_{k-1,k}$ is the correlation index
between predicted values of $\al_{it}$ under model MLAR($k-1$) and model MLAR($k$);
in boldface are the differences between maximum values of consecutive log-likelihoods
which are smaller than 0.001 for the first time.}\label{Table_choice_kq}\vspace*{0.5cm}
\end{table}

On the basis of the results in Table \ref{Table_choice_kq}, we conclude that
the adequate number of quadrature points ($q$) is
equal to $51$ for $k=1$ and to $61$ for $k=2$ and $k=3$. For
illustrative purposes, in the table we also
show results until $q = 101$. Indeed, we observe that
increasing $q$ over the selected value is unnecessary, as the
corresponding values of the maximum log-likelihood become stable. Moreover,
being $\rho_{12}$ smaller than 0.99 and $\rho_{23}$
equal to 0.9974, we choose $k = 3$ mixture components. As mentioned at
the end of Section \ref{sec:model_selection}, we note that the
proposed selection criterion for $k$ leads to
selecting a more parsimonious model with respect to that selected by
BIC. Indeed, in this last case we obtain decreasing values of the
BIC index at least until $k = 4$. Moreover, in our application this
criterion becomes soon hardly to apply, as for $k\geq 4$ the
log-likelihood becomes rather 
flat and, therefore, estimates result
highly unstable.

To evaluate the goodness-of-fit of the selected model, we compare
it with the MLAR(1) (or LAR) model and the MLAR(2) model.
We also compare these models with the LM model with covariates and initial
distribution of the Markov chain equal to the stationary distribution.
For a given number of latent states $k$, the last model is indicated by
LM($k$); it is fitted for $k=1,\ldots,10$. The results of this
comparison in terms of maximum log-likelihood and BIC index are
reported in Table \ref{tab:LM}.

\begin{table}[!ht]\centering
{\small
\begin{tabular}{l|ccc|ccc}
\hline \hline
         &   \multicolumn{3}{c|}{LM($k$)}    &   \multicolumn{3}{c}{MLAR($k$)} \\
\hline
$k$   &       \multicolumn1c{log-likel.} &  \# param. &  \multicolumn{1}{c|}{BIC} & \multicolumn1c{log-likel.} &  \# param. &  \multicolumn1c{BIC}     \\
\hline
1   &   -80,792 &   \;\;8   &   161,650 &     -63,625    &   10      &     127,340  \\
2   &   -69,866 &   \;11    &  139,830  &   -62,999 & 13 &    126,110   \\
3   &   -65,815 &   \;16    &  131,770  &   -62,846 & 16  &     125,830 \\
4   &   -64,007 &   \;23    &   128,220 &   --  &  --  &  --  \\
5   &   -63,370 &   \;32    &   127,020 &   --  &  --  &  --  \\
6   &   -63,098 &   \;43    &  126,580  &   --  &  --  &  --  \\
7   &   -63,020 &   \;56    &  126,540  &   --  &  --  &  --  \\
8   &   -62,852 &   \;71    &  126,330  &   --  &  --  &  --  \\
9   &   -62,782 &   \;88    &   126,340 &   --  &  --  &  --  \\
10  &   -62,617 &   107& 126,180    &   --  &  --  &  --  \\
\hline \hline
\end{tabular}}
\caption{\em Log-likelihood, number of parameters, and BIC index for
the LM($k$) model, with $k =1,\ldots,10$, and the MLAR($k$) model,
with $k=1,2,3$.}
\label{tab:LM}
\vspace*{0.25cm}
\end{table}

From Table \ref{tab:LM} we conclude that the smallest BIC index is for the
MLAR(3) model, to which correspond a maximum log-likelihood of -62,846 with
16 parameters. To obtain a higher log-likelihood
with the LM model, we need at least $k = 9$ latent states
and, consequently, at least 88 parameters. This confirms that
the proposed model reaches levels of goodness-of-fit comparable
with those of the LM model but, at the same time, a level of parsimony
close to that of the LAR model, since only 6 parameters are added
to this model. However, in comparing the MLAR model with
the LM model we have to consider that, especially for large values of $k$,
the likelihood of the second presents several local maxima. Therefore,
it is not ensured that the reported values of the log-likelihood for this
model corresponds to global maxima. On the other hand, at least for this application,
we did not find evidence of more local maxima of the MLAR model log-likelihood.
This is reasonable because of the reduced number of parameters.
\subsection{Parameter estimates and prediction of latent effects}
%
The estimates of the parameters of most interest of the selected
model and of comparable models are reported in Tables
\ref{tab:est_fixed} and \ref{tab:est_latent}. We recall the we
selected model MLAR(3). In particular, for models MLAR(1) (or LAR),
MLAR(2), and MLAR(3), Table \ref{tab:est_fixed} reports the
estimates of the cutpoints 
and of the regression coefficients
entering equation (\ref{eq:assumption1}) together with the
corresponding standard errors. Moreover, Table \ref{tab:est_latent}
reports the estimates of the parameters on which the latent
structure depends.

\begin{table}[!ht]\centering
{\small
\begin{tabular}{lcccc}
\hline \hline
                 & \multicolumn1c{MLAR(1)} &    \multicolumn1c{MLAR(2)}       & \multicolumn1c{MLAR(3)}  \\
\hline
$\hat{\mu}_1$   &  \;7.0155    &          \;8.2661            & \;8.9678  \\
$\hat{\mu}_2$   &  \;3.8670     &        \;4.5033             &   \;4.8351  \\
$\hat{\mu}_3$   &  \;0.6607    &         \;0.6173             &   \;0.7006  \\
$\hat{\mu}_4$   & -2.7646         &    -3.5160           &   -3.7157  \\\hline
$\hat{\be}_1$ (female)    &   -0.2056    &    -0.2148              &  -0.2317 \\
          &   (0.0738)           &     (0.0890)            &   (0.0958)    \\
$\hat{\be}_2$ non white &   -1.5175       &   -1.6884          &   -1.9735   \\
          &   (0.0968)         &    (0.1149)            &   (0.1266)    \\
$\hat{\be}_3$ education &  \;1.8182       &   \;2.2052         &   \;2.3846  \\
          &  (0.0755)          &     (0.0936)           &   (0.1020)    \\
$\hat{\be}_4$ age       &  -0.1085       &         -0.1197                &  -0.1299  \\
          &  (0.0028)          &    (0.0033)             &    (0.0036)    \\
\hline \hline
\end{tabular}
\caption{\em Fixed part estimates from  fitting 
model MLAR($k$),
$k=1,2,3$; 
in brackets are the standard errors for the regression
coefficients. \label{tab:est_fixed}}} \vspace*{0.25cm}
\end{table}

\begin{table}[!ht]\centering
{\small
\begin{tabular}{lccccc}
\hline \hline
 $k$    &   $h$    &  $\hat{\xi}_h$    &   $\hat{\rho}_h$    &   $\hat{\pi}_h$    
 &   $\hat{\sigma}^2$  \\
\hline
1       &   1    &  \;0.0000    &   0.9529    &   1.0000    &   \;9.8151  \\
\hline
2       &   1    &  -0.1073    &   0.9788    &   0.7634    &     16.5169  \\
2       &   2    &  \;0.3461    &   0.5584    &   0.2366    &     16.5169  \\
\hline
3       &   1    &  -2.6400    &   0.5204    &   0.1399     &     17.0247  \\
3       &   2    &  -0.0578    &   0.9761    &   0.7146     &     17.0247  \\
3       &   3    &  \;2.8237    &   0.3472    &   0.1455     &     17.0247  \\
\hline \hline
\end{tabular}}
\caption{\em Estimates of the parameters affecting the latent
process distribution under model MLAR($k$), $k =1, 2, 3$.
\label{tab:est_latent}}\vspace*{0.15cm}
\end{table}

From Table \ref{tab:est_fixed} we observe that the estimated
cutpoints $\hat{\mu}_1,\ldots,\hat{\mu}_4$
are ordered as we may expect in accordance with the parameterisation defined by (\ref
{eq:global}). Moreover, on the basis of the $t$-statistics
that may be computed for the regression
coefficients, we conclude that all covariates are significant under every model
considered in the table. However, the magnitude of each point estimate
increases (in absolute value) as $k$ goes from 1 to 3, while retaining the same
sign. For instance, the effect of education 
increases from 1.8182 to 2.3846.
Less evident is the variation of effect of gender (female with
respect to male), which changes from -0.2056 to -0.2317.

The results in Table \ref{tab:est_latent} imply that the estimated latent
structure is rather different under models MLAR(1),
MLAR(2), and MLAR(3). While under the MLAR(1)
model all subjects are concentrated in only one class 
characterised by a very
high correlation ($\hat{\rho}_1=0.9529$), the situation is more
complex under the other two models.  With the MLAR(2) only the
$76.34\%$ of subjects belong to a class with a high
correlation ($\hat{\rho}_1=0.9788$), whereas for the remaining $23.66\%$
the correlation estimate is located to a positive intermediate level
($\hat{\rho}_2=0.5584$), being the estimate of support point higher in
this latter class ($\hat{\xi}_2 = 0.3461$ versus $\hat{\xi}_1 = -0.1073$). Under
the MLAR(3) model, we observe one class, including the $71.46\%$
of subjects, very similar to that detected under the MLAR(1) model, being the
correlation coefficient equal to 0.9761 and the corresponding
support point very close to 0. Then, the remaining part of subjects
results equally distributed between the other two classes, that show
opposite values of the support points ($\hat{\xi}_1 = -2.6400$ and $\hat{\xi}_3 =
2.8237$) and intermediate levels of correlation. Mostly
for subjects in class 3  the correlation between individual effects in
consecutive occasions is rather weak ($\hat{\rho}_3=0.3472$), so that we
may suppose that these subjects are characterised by sudden
changes in unobservable factors affecting their health status.

\begin{figure}[!h]\centering
\begin{tabular}{cc}
\includegraphics[width=8cm]{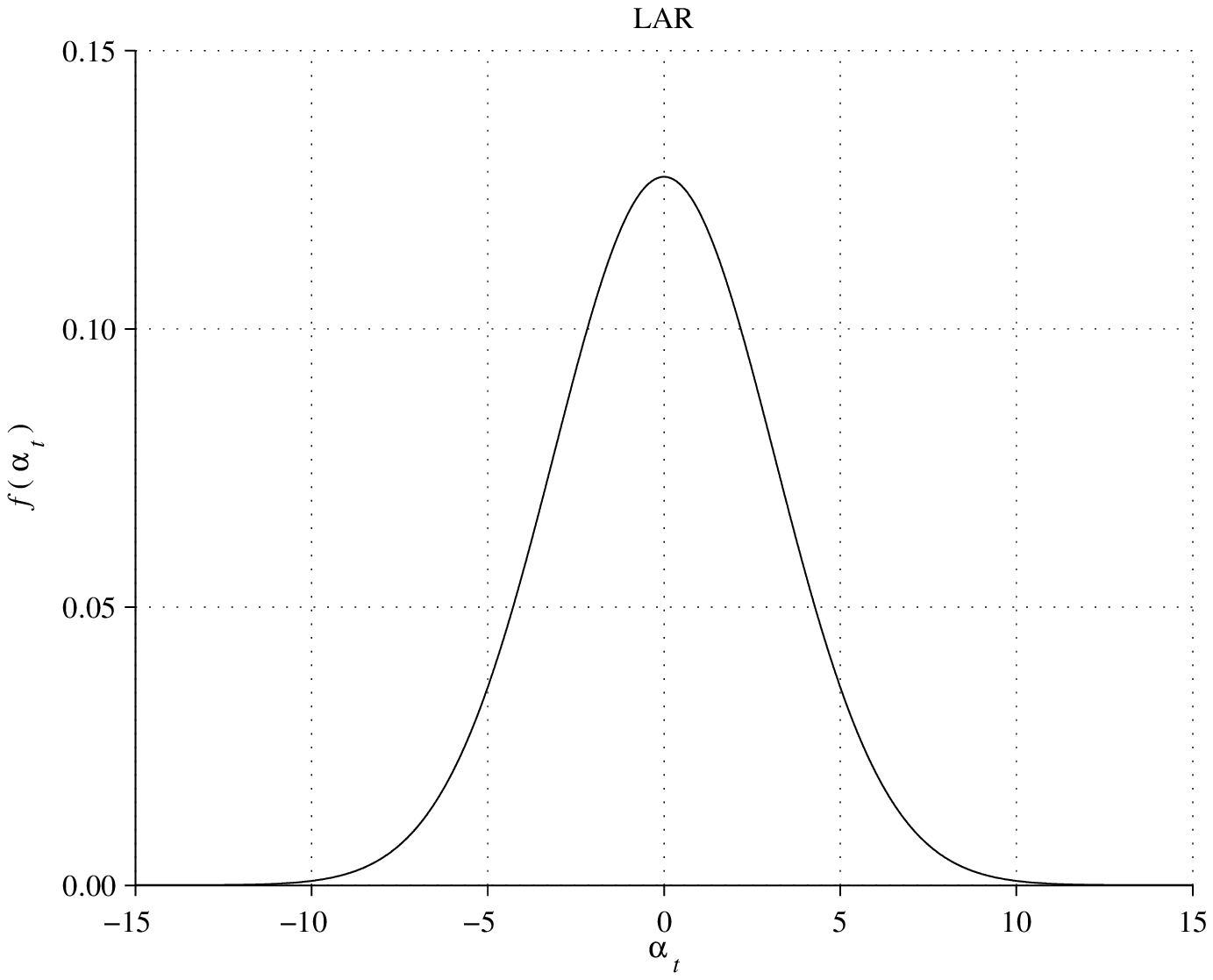} &
\includegraphics[width=8cm]{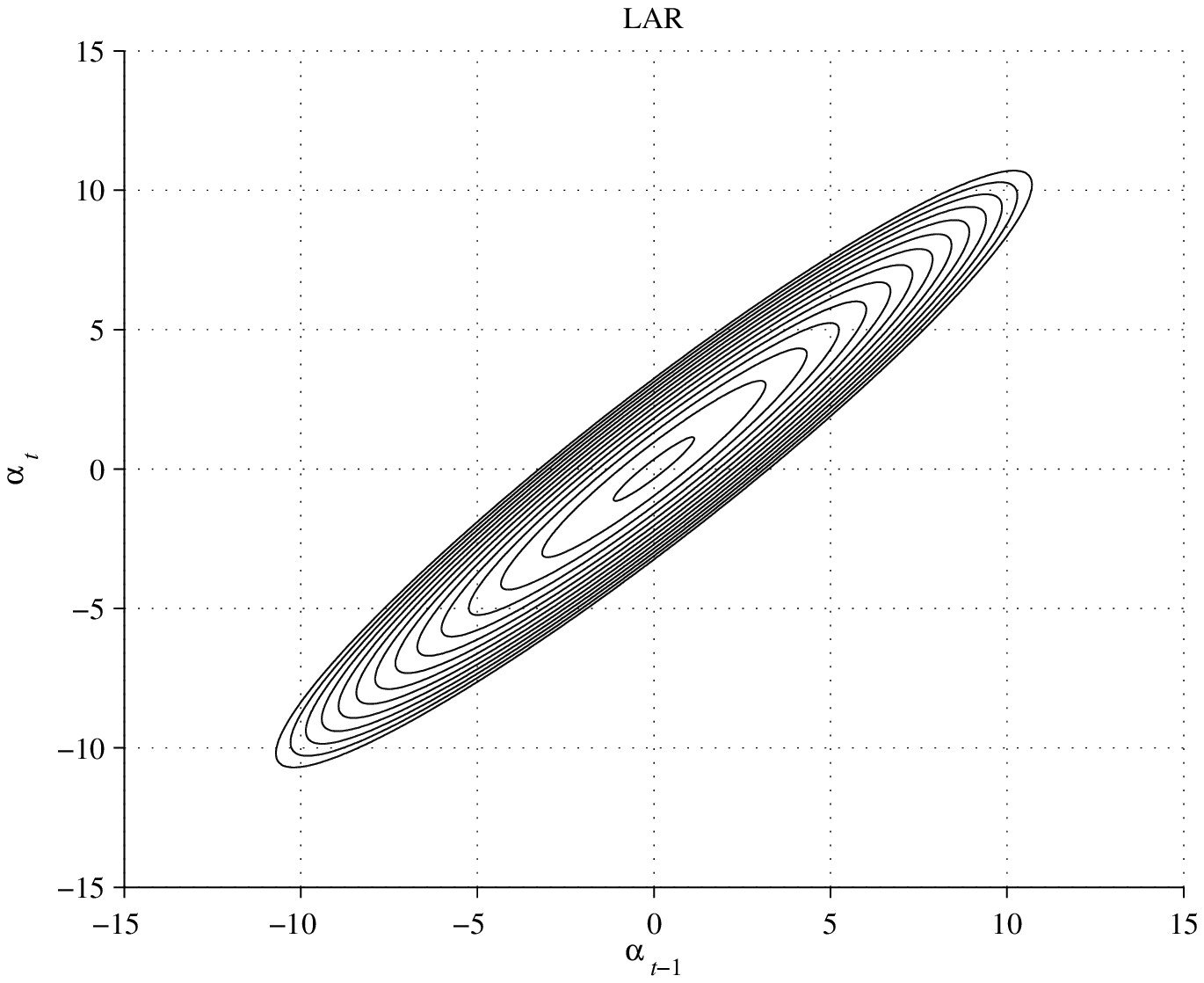}\\
\includegraphics[width=8cm]{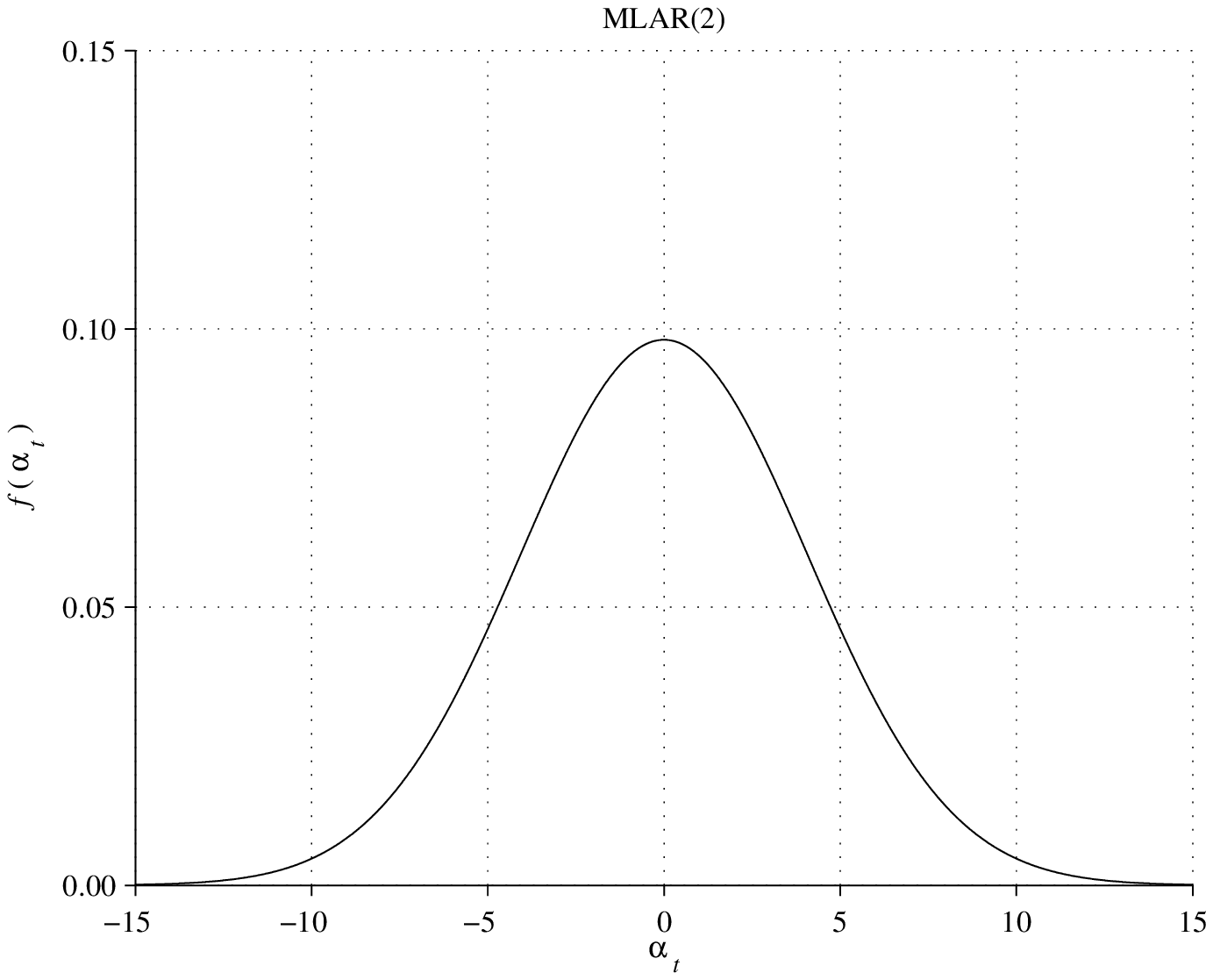} &
\includegraphics[width=8cm]{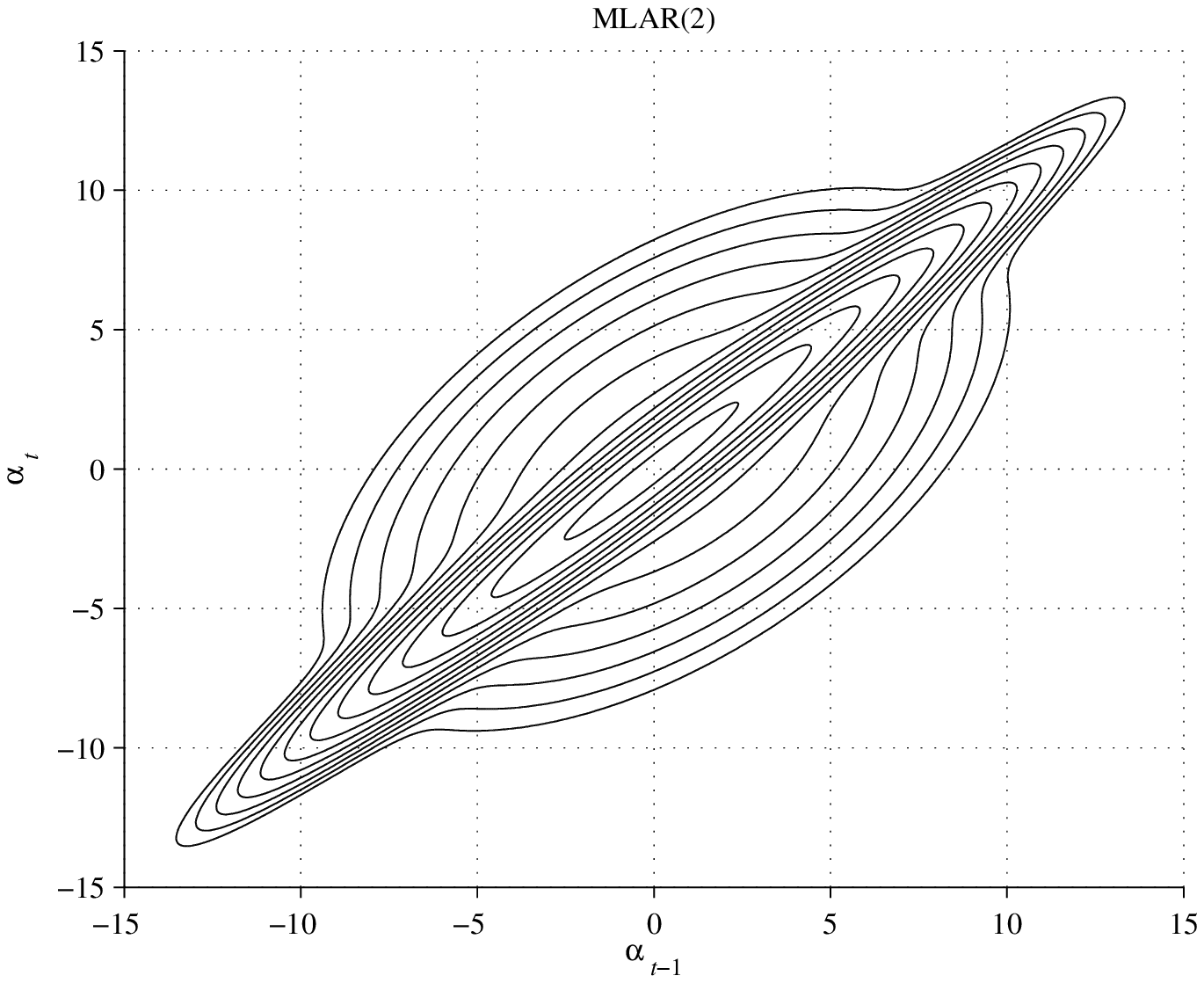}\\
\includegraphics[width=8cm]{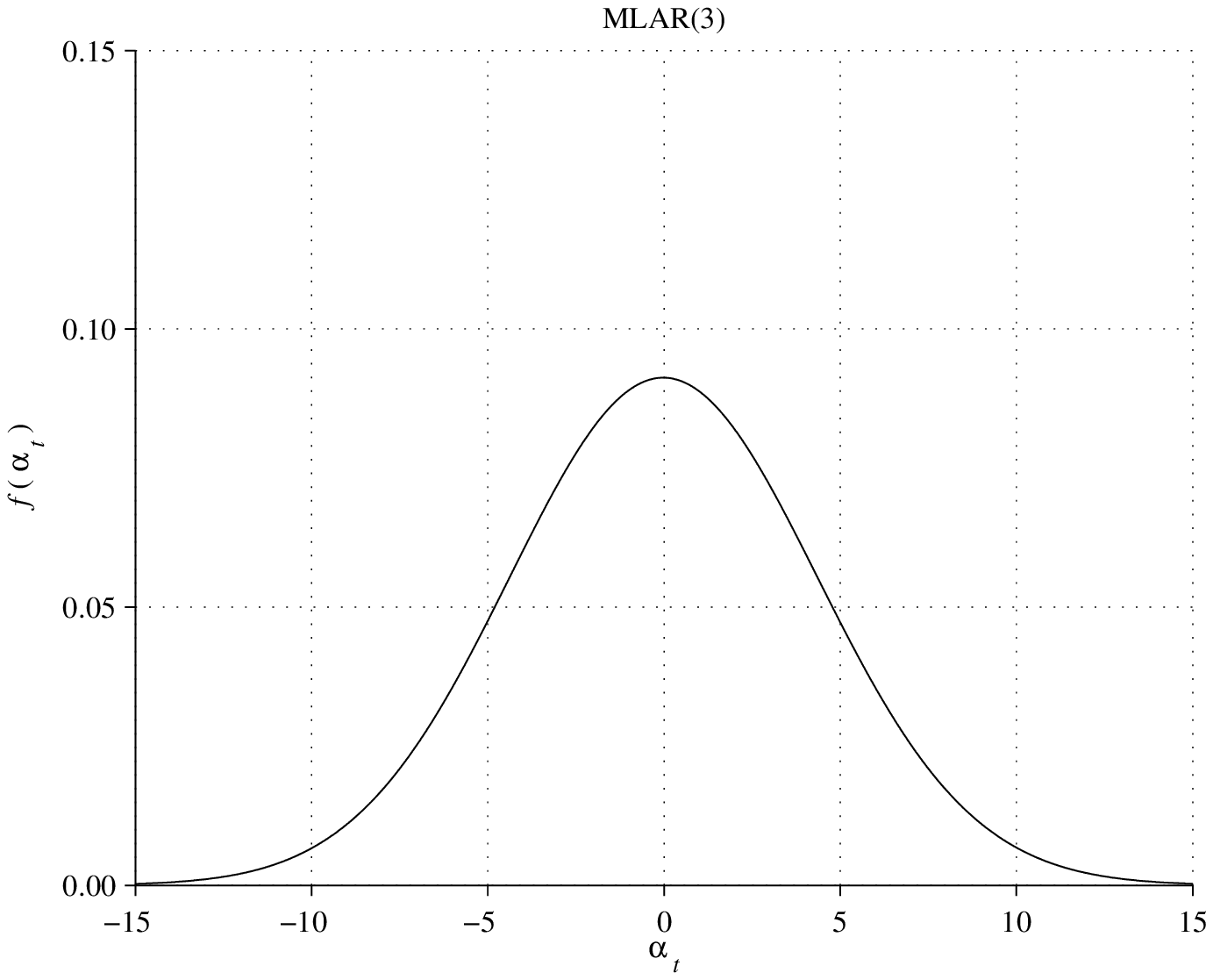} &
\includegraphics[width=8cm]{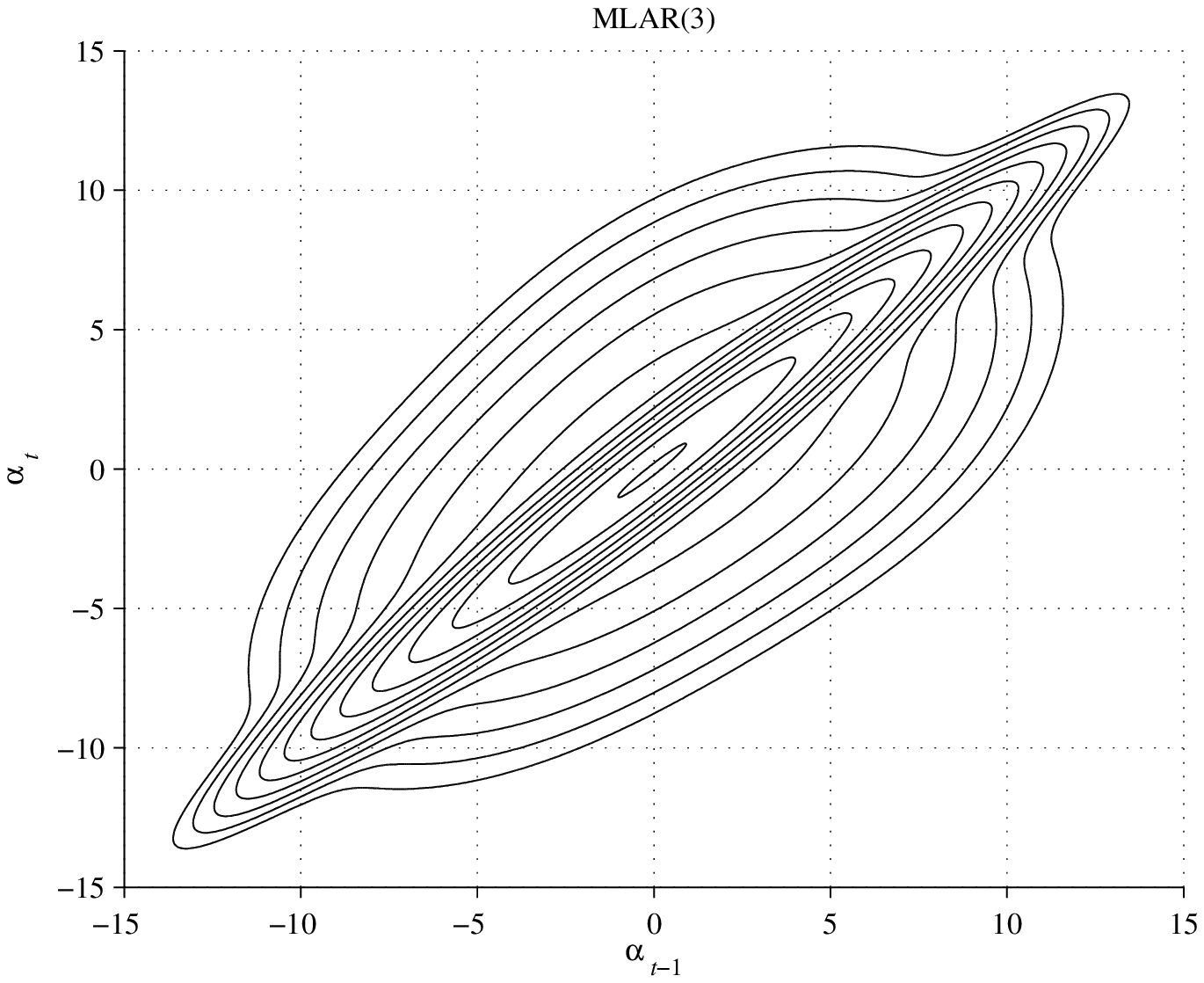}
\end{tabular}
\caption{\em Estimated density function of the univariate distribution of $\al_{it}$
(left) and of the bivariate distribution of $(\al_{i,t-1},\al_{it})$
(right) under the LAR model (top), the MLAR(2) model (middle), and the
MLAR(3) model (bottom).}
\label{fig:appl}
\end{figure}

The above results are well illustrated
by the density functions of the univariate distribution of
the individual effects $\al_{it}$ and of the bivariate distribution
of $(\al_{i,t-1},\al_{it})$ represented in Figure
\ref{fig:appl}. Indeed, while in the MLAR(1) model values at time
$t$ are highly correlated to values at time $t-1$, under the MLAR(2) model
and, more evidently under MLAR(3) model, a higher dispersion of the individual
effects $\al_{it}$ is observed on the bivariate plot. On the other hand,
the univariate distributions does not seem to deviate from normality
under both MLAR(2) and MLAR(3) models.

\begin{figure}[!h]\centering
\begin{tabular}{cc}
\includegraphics[width=8cm]{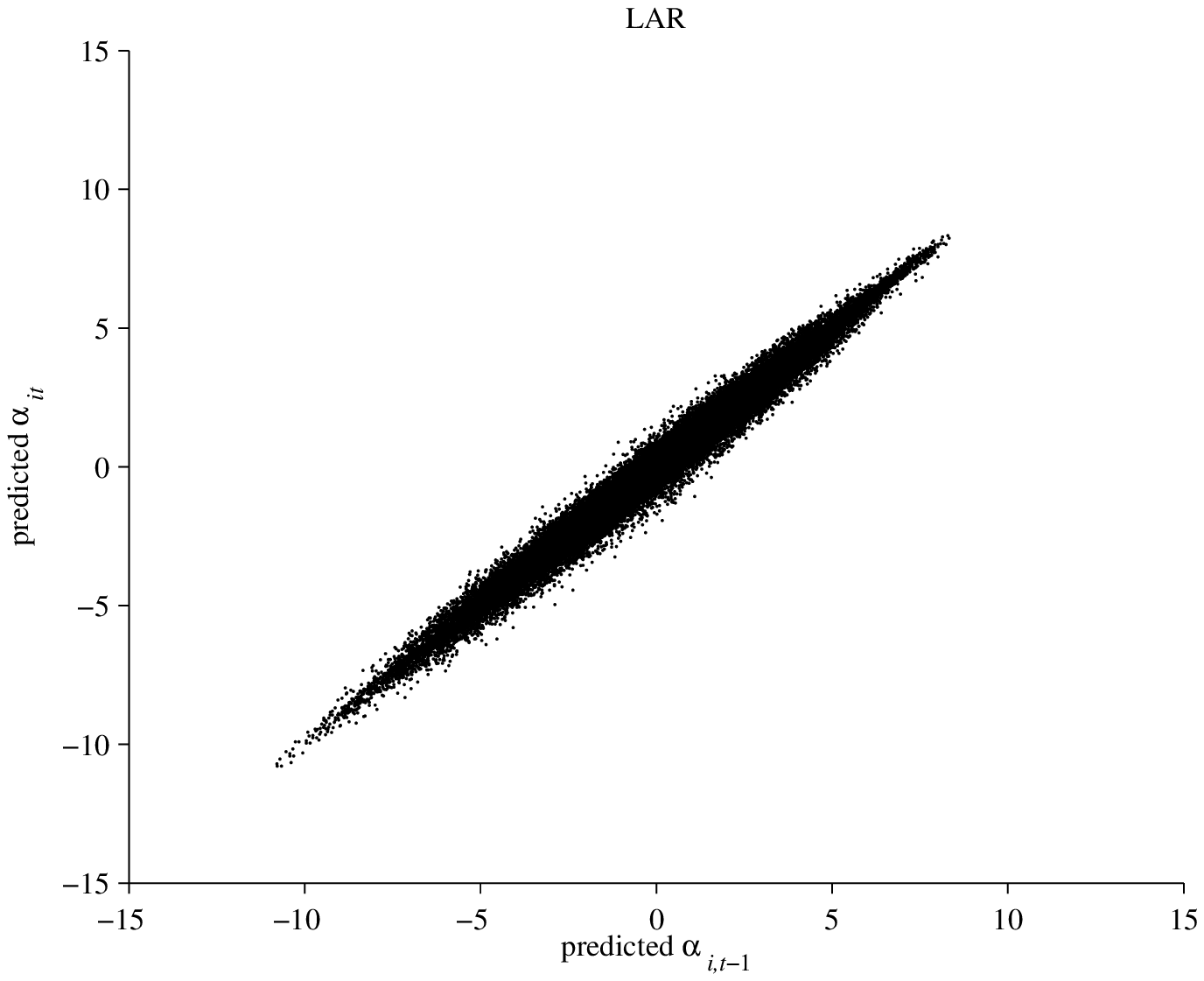} &
 \includegraphics[width=8cm]{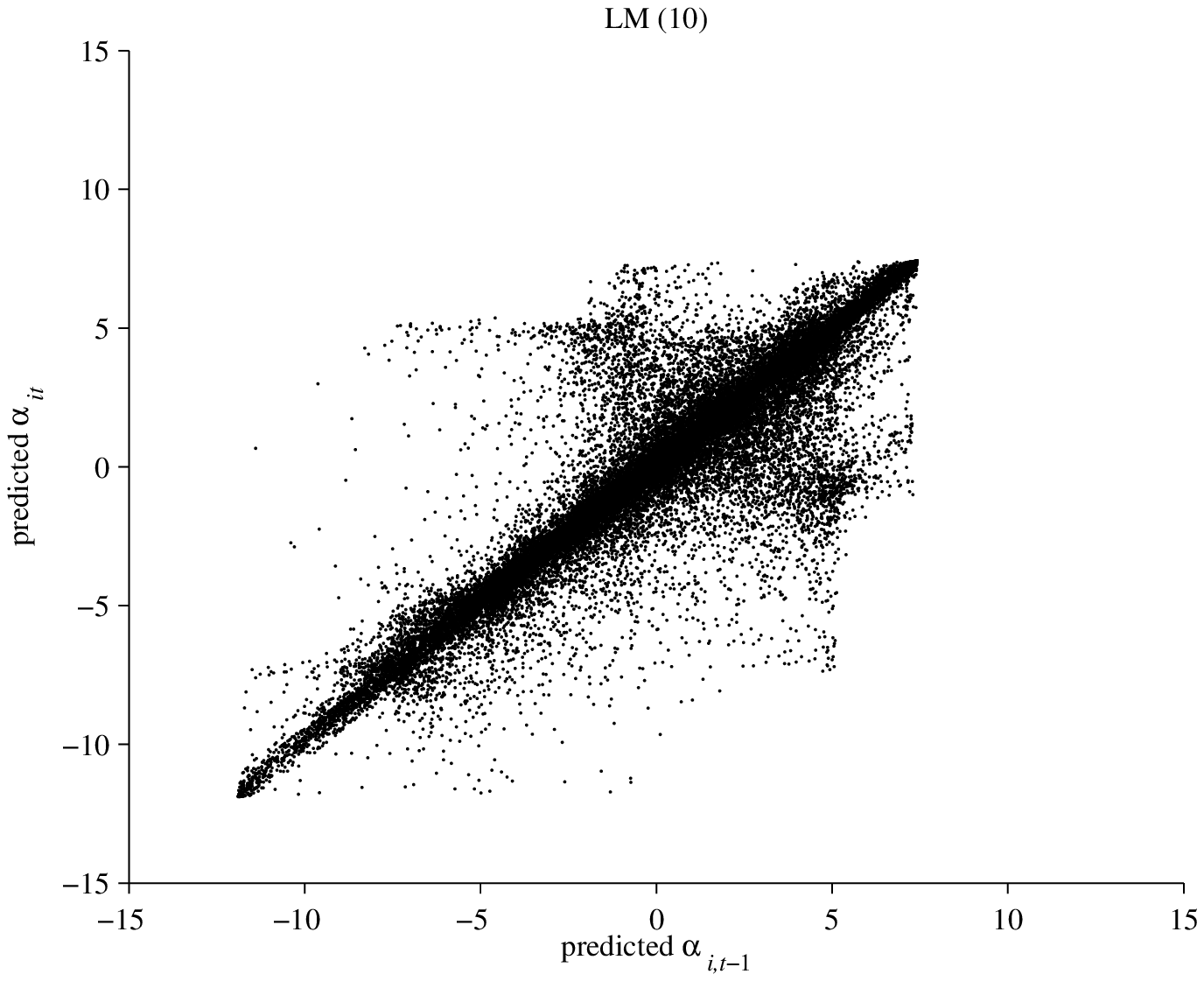} \\
 \includegraphics[width=8cm]{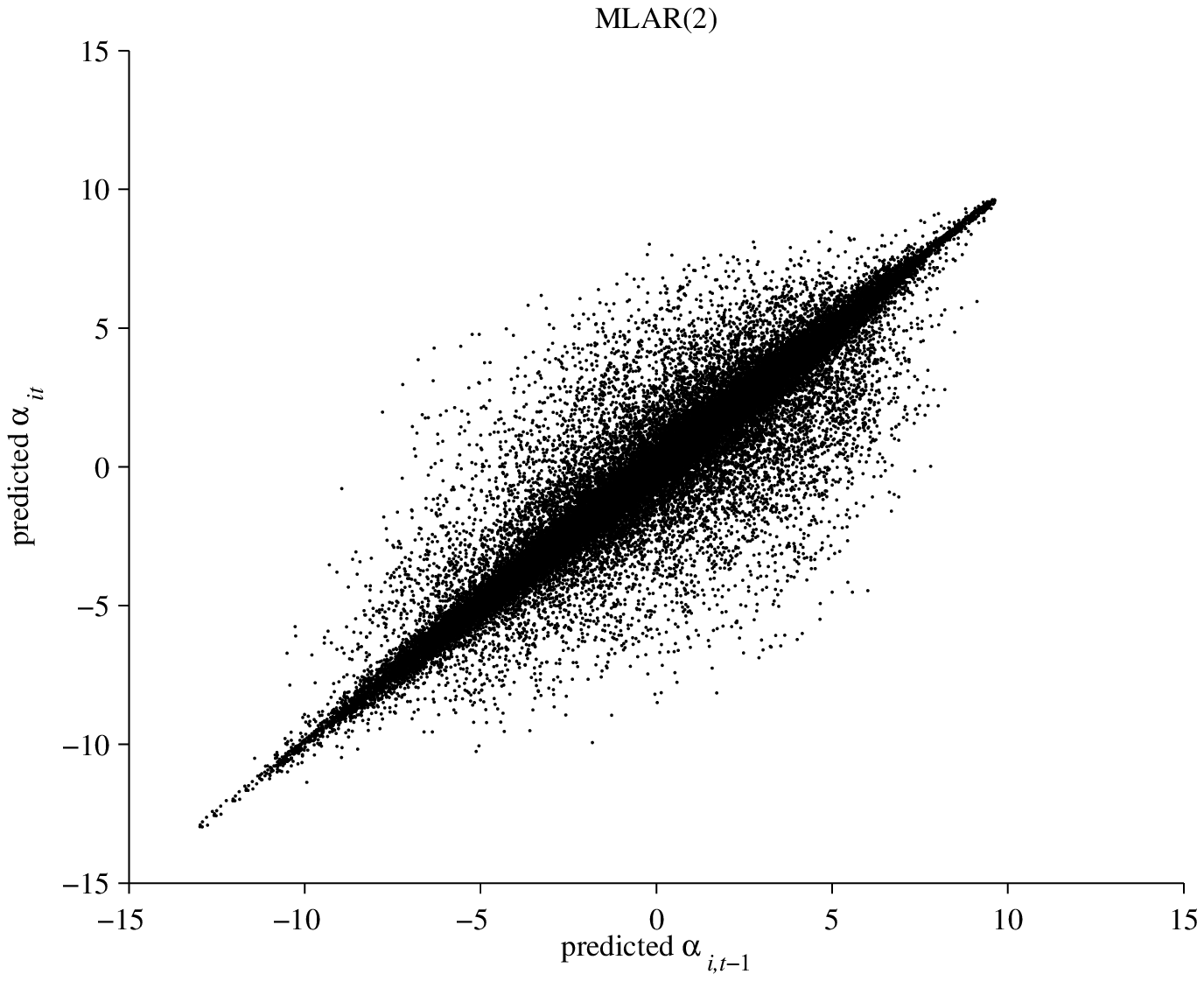} &
 \includegraphics[width=8cm]{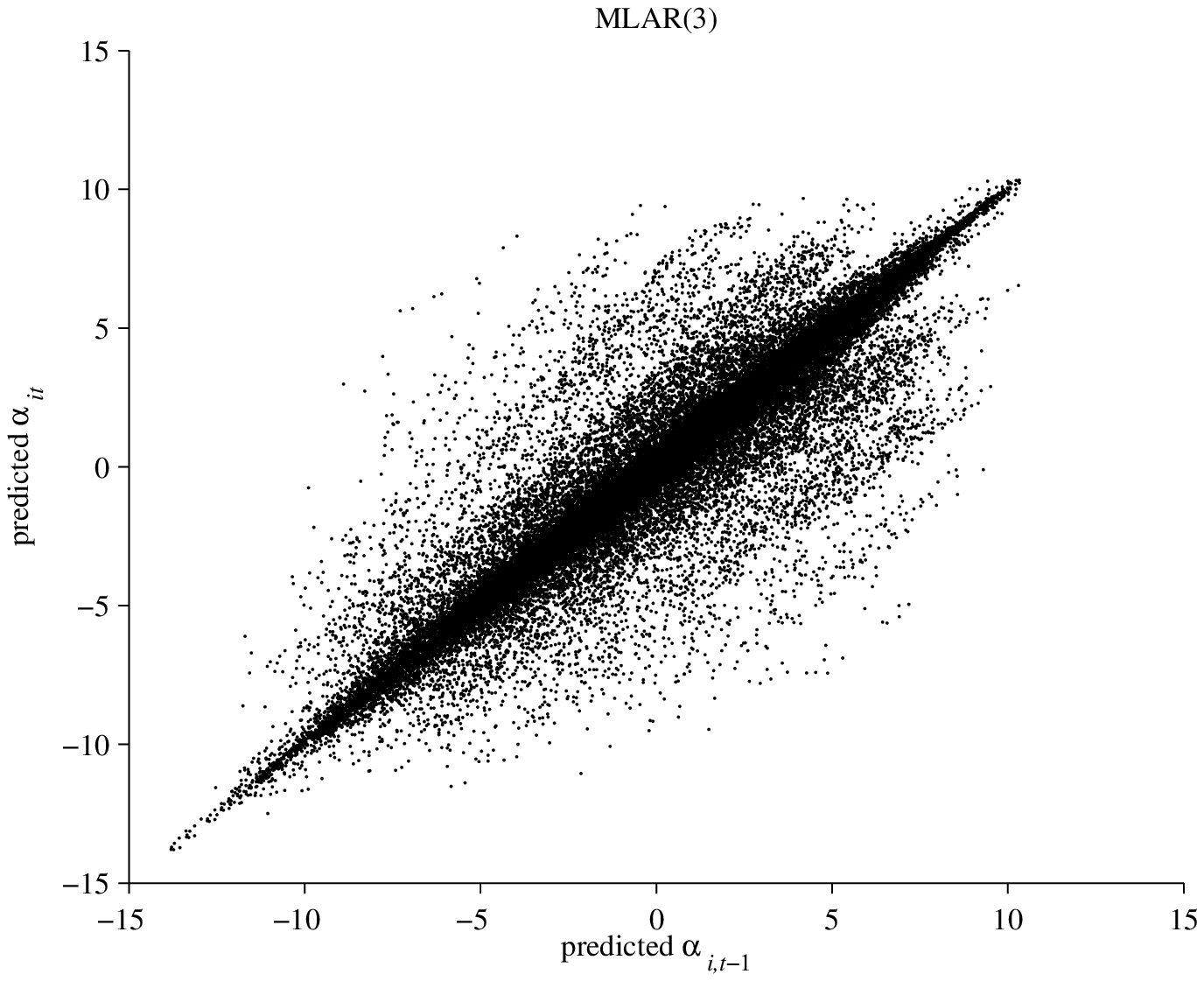}
\end{tabular}
\caption{\em Predicted values of $(\al_{i,t-1},\al_{it})$, $i=1,\ldots,n$,
$t=2,\ldots,T$, under the LAR model (top  left), the LM model (top right),
the MLAR(2) model (bottom left), and the MLAR(3) model (bottom right).}
\label{fig:predict}
\vspace*{0.5cm}
\end{figure}

Another way to compare the estimated latent structure under the different
models here considered is through the predicted $\al_{it}$, as
computed by (\ref{eq:pred}), rather than
through the a priori distributions in the previous Figure. These predicted
values are represented, also for the same LM(10) model already considered in
Section \ref{sec:application:model_selection}, in Figure \ref{fig:predict}.
In particular, each plot in the figure represents points with coordinates
($\hat{\al}_{i,t-1},\hat{\al}_{it}$), for $i=1,\ldots,n$ and $t=2,\ldots,T$.

According to the plots in Figure \ref{fig:predict}, under models MLAR(2), MLAR(3),
and LM(10) there is a stronger dispersion of the predicted
values with respect to the LAR model. In particular, the difference of the
plot obtained under the selected MLAR(3) model is neat with respect to the
plot obtained under the LAR model, whereas it is less evident with respect to the
plot obtained under the LM(10) model, which however uses much more parameters.
In practice, it seems that the proposed
model allows for more erratic trends of the unobserved
individual effects across time with
respect to the LAR model. This is a direct consequence of the greater flexibility
of the MLAR model, which is obtained at the cost of
a reduced number of additive parameters.\vspace*{-0.25cm}
\section{Discussion}
In this paper, we extend the Latent Autoregressive (LAR) model for
longitudinal data  \citep{chi:89, hei:08}, by adopting a latent
structure for the unobserved heterogeneity which is based on a
mixture of AR(1) processes with specific mean values and correlation
coefficients, but with common variance. The proposed model, named
Mixture Latent Autoregressive (MLAR) model, is formulated in a
general way, so that it may be easily adapted to different types of
response variable (binary, ordinal, or continuous). It is important
to note that, as for the LAR model, the latent process on which the
proposed model is based is continuous.

Compared to the latent
Markov (LM) model with covariates in which the latent process is discrete
(see \cite{bart:farc:penn:10} for a review),
the MLAR model has some interpretative advantages, being usually more
natural to consider the effect of unobservable factors or covariates
as continuous rather than discrete. Moreover, the main advantage of MLAR
model with respect to LAR model is the improvement of the goodness-of-fit
which becomes close to that of an LM model with covariates,
allowing us to adequately take
into account more erratic trends in the unobservable individual effects.
At the same time, the parsimony of the proposed model
is kept near to that of a LAR model, avoiding the explosion of the
number of parameters of the LM model when the number of latent
states increases.

In order to make inference on the proposed model, we show how its
likelihood may be efficiently computed by exploiting some
recursions developed in the context of Hidden Markov (HM) models
\citep{bau:70,mac:zuc:97}. The resulting computational method is
equivalent to the Sequential Gaussian Quadrature method proposed by
\cite{hei:08} for the LAR model. 
Moreover, since the model likelihood is the same as that of an LM model
with suitable constraints, then maximum likelihood estimation of the model parameters 
may be performed on the basis of an EM algorithm for LM models, adapting that 
implemented by
\cite{bart:farc:09}. 
To make faster the estimation, after a
certain number of EM steps we suggest to switch to a Newton-Raphson algorithm, which
is based on the observed information matrix obtained by a numerical method.

After parameters estimation, standard errors are obtained from the
observed information matrix. Moreover, on the basis of the parameter
estimates it is also possible to predict individual effects for
every subject and time occasion. We show how these predicted values
may be used to implement a model selection strategy for the number
of mixture components of the MLAR model; we recall that each
component corresponds to a separate AR(1) latent process. In
particular, the number of mixture components is increased until the
predicted values of the latent variable do not significantly change.
This selection strategy leads us
to selecting more parsimonious models with respect to alternative
methods, such as those based on information criteria.

The advantages of the MLAR model with respect to the LAR model are
illustrated through an application to a longitudinal dataset, coming
from the Health and Retirement Study conducted by the University of
Michigan, about self-evaluation of the health status. The results
show evidence of three mixture components corresponding to the same
number of AR(1) processes. 
Each component has its own specific correlation parameter. 
In this way, we take into account the latent trends and
jumps from time to time in a more flexible way in comparison with
the LAR approach. Moreover, the goodness-of-fit 
of the MLAR model is comparable to
that of an LM model with covariates, being, at the same time, the
number of parameters strongly smaller. In addition, we observe that the proposed model 
does not suffer from
the problem of multimodal likelihood as the LM model with covariates
does. This is rather obvious considering the reduced number of
parameters of the first with respect to the second.

A final point concerns possible extensions of the proposed approach.
A natural extension consists of generalising our model to AR processes
of order two (or higher), so as to take also into account 
more sophisticated dependence structures.
This extension may be performed along the same lines as in \cite{bart:solis:10},
although we think that likelihood inference
becomes more problematic, especially in terms of computational capability
required to obtain maximum likelihood estimates.

Finally, the reader may wonder about a possible extension of the
proposed approach in which subjects may move between different
mixture components or, in equivalent terms, in which we have
switching parameters for the individual AR(1) latent processes. This
amounts to combine a latent AR(1) process with a latent Markov chain.
We recall that in the MLAR model the latent process 
parameters are kept 
constant across time, although these parameters may be different among subjects. The
extended formulation based on combining a latent AR(1) process and a
latent Markov chain has been experimented by the same authors, see
\cite{bacc:bart:penn:10}; however, they observed that the major
complexity of the resulting model tends to produce estimates that
are highly unstable and are rather difficult to interpret. Moreover,
the goodness-of-fit reached under this extended formulation is
comparable to that of the 
MLAR model here proposed with the same number of parameters. For
these reasons we consider this model as the right compromise between
making more flexible the LAR model and keeping a parsimonious
structure, while retaining a continuous latent process approach.
\section*{Acknowledgments}
F.Bartolucci tanks Prof. F.Peracchi for stimulation discussions on
the topic and, together with F.Pennoni, acknowledge the financial
support from the ``Einaudi Institute for Economics and Finance''
(EIEF), Rome (IT).
\bibliography{biblio}
\bibliographystyle{apalike}

\end{document}